\documentclass{amsart}
\usepackage[dvips]{graphics}

\newtheorem*{proposition*}{Proposition}
\newtheorem*{theorem*}{Theorem}
\newtheorem{theorem}{Theorem}[section]

\newtheorem{proposition}[theorem]{Proposition}
\newtheorem{lemma}[theorem]{Lemma}

\newcommand{\Diagram}{\mathcal{D}}

\newcommand{\Twist}{\tau}

\newcommand{\Crossing}{\mathrm{cr}}
\newcommand{\Diam}{\mathrm{Diam}}
\newcommand{\Edgepath}{\gamma}
\newcommand{\EdgepathSystem}{\Gamma}
\newcommand{\BasicEdgepath}{\lambda}
\newcommand{\BasicEdgepathSystem}{\Lambda}

\newcommand{\dec}{\mathrm{dec}}
\newcommand{\inc}{\mathrm{inc}}

\newcommand{\Seifert}{\mathrm{Seifert}}

\newcommand{\NumTangles}{N}

\newcommand{\EDGE}{\,--\,}
\newcommand{\DIREDGE}{\,--\,}

\newcommand{\Radius}{\rho}

\newcommand{\angleb}[1]{\langle #1 \rangle} % angle bracket

\newcommand{\circleb}[1]{\langle #1 \rangle^{\circ}} % angle bracket

\begin{document}

\title[Crossing number and diameter of boundary slope set]{Crossing number and diameter of boundary slope set of Montesinos knot}

\author{Kazuhiro Ichihara}
\address{%
College of General Education, 
Osaka Sangyo University, 
3--1--1 Nakagaito, Daito, Osaka 574--8530, Japan }
%\email{ichihara@las.osaka-sandai.ac.jp}
\curraddr{
  Department of Mathematics Education,
  Nara University of Education,
  Takabatake-cho, Nara, 630-8528, Japan
}
\email{ichihara@nara-edu.ac.jp}

\author{%
    Shigeru Mizushima}
\address{%
        Department of Mathematical and Computing Sciences \\
        Tokyo Institute of Technology \\
        2--12--1 Ohokayama, Meguro \\
        Tokyo 152--8552, Japan}
\email{mizusima@is.titech.ac.jp}

\thanks{The first author is supported in part by Grant-in-Aid for Young Scientists (B) \ (No.\ 18740038).}
% 57M25 : Knots and links in S^3

%\date{\today}

%%%%%%%%%%%%%%%%%%%%%%%%%%%%%%%%%%%%%%%%%%%
%
% Abstract
%
%%%%%%%%%%%%%%%%%%%%%%%%%%%%%%%%%%%%%%%%%%%

\begin{abstract}
It is shown that 
the diameter of the boundary slope set 
is bounded from above by 
the twice of the minimal crossing number 
for a Montesinos knot. 
\end{abstract}

\maketitle

\section{Introduction}

For a knot $K$ in the $3$-sphere $S^3$, 
the minimum number of crossings of a diagram 
among all diagrams of $K$ 
is called the \textit{crossing number}.
This is one of the most basic topological complexities of a knot. 

Another complexity of knots which we consider in this paper is 
the diameter of the set of boundary slopes. 
Let $E(K)$ be the exterior of $K$. 
A compact connected surface properly embedded in $E(K)$ 
is said to be \textit{essential} 
if the surface is incompressible and boundary-incompressible.
The boundary of an essential surface appears on $\partial E(K)$ 
as a parallel family of non-trivial simple closed curves, 
and so, it determines an isotopy class of the curves, 
which is called the \textit{boundary slope} of the surface. 
Recall that 
the set of slopes is identified with the set of rational numbers with $\infty$, 
where the meridian of $K$ corresponds to $\infty$, 
by using the standard meridian-longitude system for $K$. 
Then \textit{the diameter of the set of boundary slopes for $K$} 
is defined as 
the difference between the greatest one and the least one as rational numbers, 
except for the infinity (i.e., meridional) slope. 
This is well-defined since 
it is known that there are only finitely many boundary slopes \cite{H} 
and there exist at least two boundary slopes in general \cite{CS84}. 

Our main result is the following: 

%
% Theorem
%
\begin{theorem}[Main Theorem]
\label{Thm:CrossingNumberAndDiameter:Main}
Let $K$ be a Montesinos knot without $1/0$-tangles.
Then, we have
\begin{eqnarray*}
2\,\Crossing(K)&\ge&\Diam(K)
,
\label{Eq:Diam:Crossing:Main}
\end{eqnarray*}
where $\Crossing(K)$ denotes the crossing number of $K$ 
and $\Diam(K)$ denotes the diameter of the set of boundary slopes for $K$.
%
%%%The equality holds for all alternating Montesinos knots. 
The equality holds if $K$ is alternating. 
\end{theorem}

It is probable that this inequality holds for all knots in general. 
For example, it is easily verified for torus knots. 
For a non-trivial torus knot $T_{p, q}$, 
assuming that $p, q$ are relatively prime with $2 \leq q \leq p$, 
it is known that the boundary slopes are $0$ and $pq$ 
while the crossing number is $pq- p$. 
Thus, we have 
$2\,\Crossing(T_{p, q}) = pq + p(q - 2) \geq pq = \Diam(T_{p, q})$. 
The equality holds when $q=2$, equivalently, the knot is alternating. 
Moreover, as far as the authors observed, 
the equality holds for all alternating knots. 
We here remark that, for alternating knots, 
the converse inequality $\Diam(K) \geq 2\,\Crossing(K)$ can be shown. 
By the checkerboard construction for the reduced alternating diagram,
we can obtain two surfaces, which are both essential by \cite{A,DR}, 
and the difference between their boundary slopes is actually $2\,\Crossing(K)$. 

To prove the theorem, 
we mainly use the algorithm given by Hatcher and Oertel  in \cite{HO}, 
which enumerates all boundary slopes for a given Montesinos knot. 
The algorithm in turn is based on 
the algorithm by Hatcher and Thurston in \cite{HT},
which enumerates all boundary slopes for a two-bridge knot.
Another ingredient is the way to calculate 
the crossing number of a given Montesinos knot 
obtained by Lickorish and Thistlethwaite in \cite{LT}. 
In the next section, 
we introduce some notation and terminology.
Then, we prove the main theorem in the last section.

The authors would like to thank Masaharu Ishikawa 
for suggesting the potential relationship between the crossing number and 
the diameter of the set of boundary slopes. 
They also thank to Thomas Mattman for 
letting them know the work \cite{MMR}. 
They are grateful to the referee for careful reading of the manuscript and useful comments

%%%%%%%%%%%%%%%%%%%%%%%%%%%%%%%%%%%%%%%%%%%
%
% Section : Preliminary
%
%%%%%%%%%%%%%%%%%%%%%%%%%%%%%%%%%%%%%%%%%%%

\section{Preliminary}

In this section, we prepare various notations and recall fundamentals.

\begin{figure}[hbt]
 \begin{picture}(140,65)
%1over2-tangle-on-sphere-copy.eps[14] PS 78x78 78x78+0+0
%1over2-tangle.eps[16] PS 104x104 104x104+0+0
  \put(0,8){\scalebox{0.6}{\includegraphics{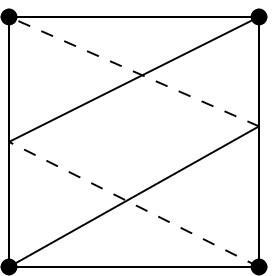}}}
  \put(80,0){\scalebox{0.6}{\includegraphics{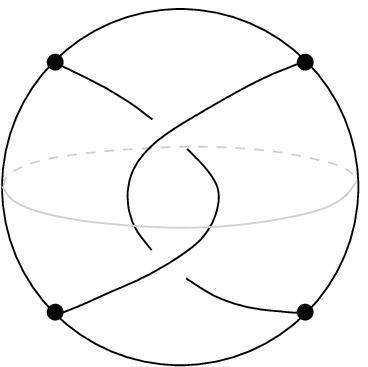}}}
 \end{picture}
 \caption{Arcs on a pillowcase and $1/2$-tangle}
 \label{Fig:1over2-Tangle}
\end{figure}

We start with a {\em rational tangle}, which we regard as follows. 
See Figure \ref{Fig:1over2-Tangle}.
The left figure illustrates a sphere with four punctures, which is flattened like a ``pillowcase''.
%%% In the left figure,
Consecutive segments with slope $1/2$ are drawn from each of four punctures.
These segments form two arcs.
Think about the interior of the sphere,
and push these two arcs into the interior with four ends fixed.
See the right figure in Figure \ref{Fig:1over2-Tangle}.
Then, we have a $1/2$-tangle in a $3$-ball.
Similarly, we can obtain a $p/q$-tangle 
for arbitrary irreducible fraction $p/q$ including $1/0$.

\begin{figure}[htb]
 \begin{center}
%0 0 434 362
%  \begin{picture}(130,109)
%\scalebox{0.3}
  \begin{picture}(109,91)
   \put(0,0){\scalebox{0.25}{\includegraphics{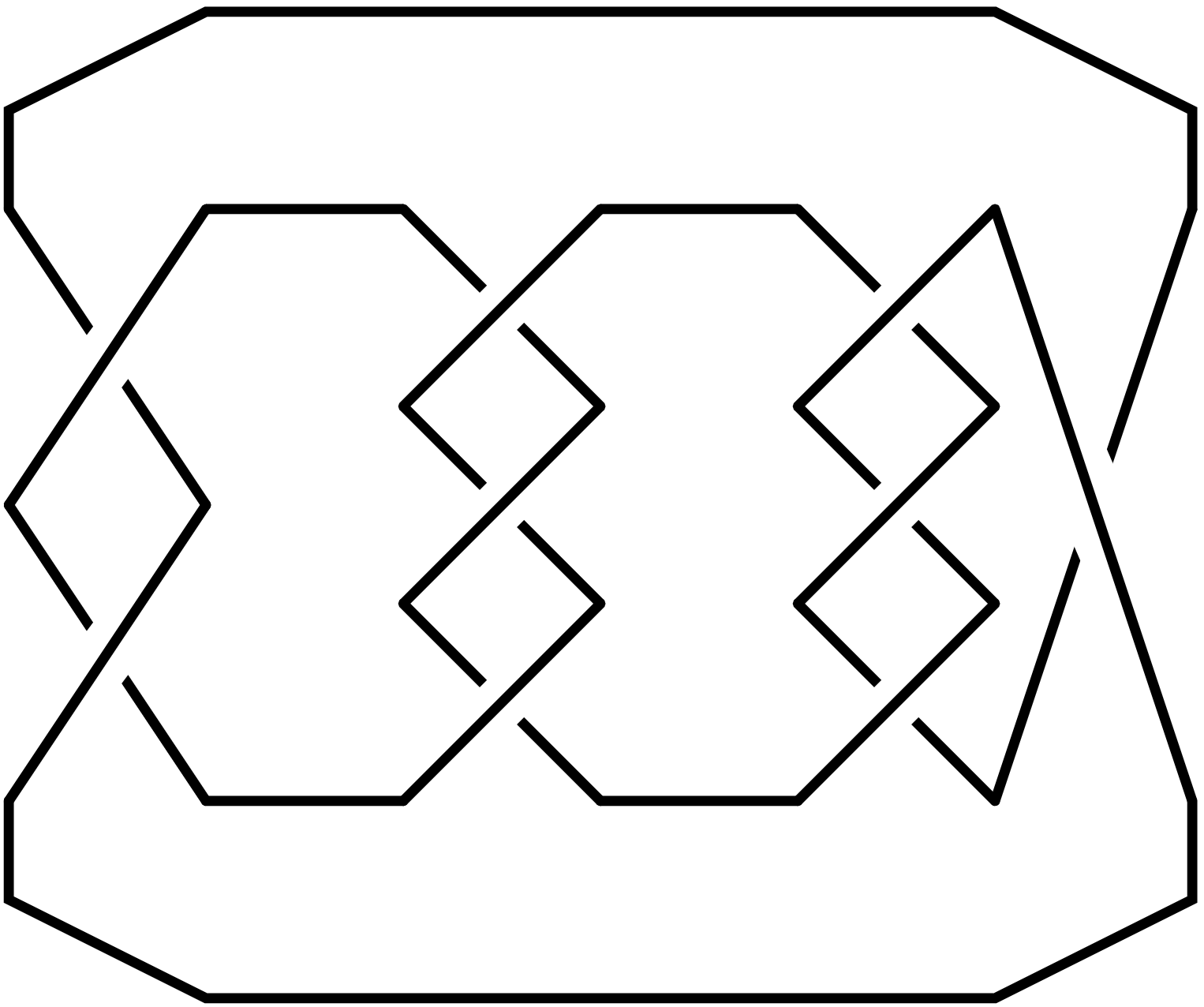}}}
  \end{picture}
  \caption{A diagram of a Montesinos knot $K(1/2,1/3,-2/3)$}
%-Montesinos knot
  \label{Fig:MontesinosKnot}
 \end{center}
\end{figure}

A {\em Montesinos knot} is defined as a knot obtained
by putting rational tangles together in a circle.
See Figure \ref{Fig:MontesinosKnot}.
A Montesinos knot obtained from
rational tangles 
$R_1$, $R_2$, $\ldots$, $R_\NumTangles$ is denoted by
$K(R_1$, $R_2$, $\ldots$, $R_\NumTangles)$.
The number of tangles will be denoted by $N$ throughout this paper.
Concerning such a tuple $(R_1$, $R_2$, $\ldots$, $R_\NumTangles)$, 
we will keep the following assumptions: 
\begin{itemize}
\item
In general, by combining rational tangles, 
we have a Montesinos link with one or more link components. 
However, in this paper, we want to only consider Montesinos knots; the number of components is one. 
Thus we assume that the tuple must satisfy either of the two conditions; 
(a) exactly one of the denominators of $R_i$'s is even
or
(b) all denominators of $R_i$'s are odd and the number of odd numerators is odd.
\item
We always assume that none of $R_i$'s are $1/0$.
%We call a Montesinos knot satisfying this assumption a {\em Montesinos knot without $1/0$-tangles}.
Furthermore, we always assume $N\ge 3$ since if the number of tangles is two or less,
the knot is found to be a two-bridge knot.
For the two-bridge knots, a result corresponding to the main theorem in this paper 
has been presented by Mattman, Maybrun and Robinson in \cite{MMR}.
\item
We assume that each tangle is non-integral, 
since an integral tangle can be combined with an adjacent rational tangle. 
This operation is a kind of normalization.
\end{itemize}

%
% Subsection
%

\subsection{Boundary slopes for Montesinos knot}

In this subsection, we review the notions developed 
in the work of Hatcher and Oertel \cite{HO} 
about the boundary slopes for Montesinos knots. 
Though we will try to give necessary explanations 
so that this paper can be self-contained, 
it would be more preferable that 
the reader is rather familiar to their machinery. 
Also see our previous note \cite{IM} for detail.

Throughout the following, 
let $K$ be a Montesinos knot obtained from rational tangles 
$R_1$, $R_2$, $\ldots$, $R_\NumTangles$. 
In sequel, $R_i$ will denote 
a fraction or the corresponding rational tangle 
depending on the situation.

\medskip

We first give an outline of the machinery developed by Hatcher and Oertel briefly. 

Take a simple unknotted loop in the 3-sphere $S^3$,
and set disjoint $\NumTangles$ disks bounded by the loop. 
By the disks, $S^3$ is divided into $\NumTangles$ 3-balls.
After appropriate isotopies, $K$ is divided into rational tangles 
%%% $(R_i)_{i=1}^{\NumTangles}$,
$(R_1, R_2, \ldots, R_\NumTangles)$,
the intersections with the $\NumTangles$ balls.
Let $F$ be an essential (meaning that, incompressible and boundary-incompressible) surface 
properly embedded in the exterior of $K$. 
By virtue of Proposition 1.1 in \cite{HO}, 
$F$ can be isotoped so that it is divided into subsurfaces 
%%% $(F_i)_{i=1}^{\NumTangles}$,
$(F_1, F_2, \ldots, F_\NumTangles)$ in certain ``good position'' 
included in each of the $\NumTangles$ balls.
Then each $F_i$ is represented by 
an ``edgepath'' $\Edgepath_i$ in certain ``diagram'' $\Diagram$, 
and the whole $F$ is represented 
by an ``edgepath system'' $\EdgepathSystem$ in $\Diagram$. 
Actually the boundary slope of $F$ is 
calculated from the corresponding $\EdgepathSystem$ 
in purely combinatorial way. 
Conversely an ``edgepath system'' $\EdgepathSystem$ in $\Diagram$ 
corresponds to a properly embedded surface in the exterior of $K$. 
In \cite{HO}, conditions for determining essentiality 
of the surface are fully described. 

\medskip

In the following five subsubsections, 
we will summarize the features of the diagram $\Diagram$, 
edges, edgepaths and edgepath systems in $\Diagram$, 
which we will use in the rest of paper.

%%%%%%%%%%%%%%%%%%%%%%%%%%%%%%%

%
%	Diagram
%

\subsubsection{Diagram}

\begin{figure}[htb]
%   \begin{center}
%    \begin{minipage}{200pt}
     \begin{center}
      \begin{picture}(70,210)
%       \put(0,0){\scalebox{0.7}{\includegraphics{diagram.eps}}}
%       \put(-14,2){\scalebox{1.0}{\includegraphics{diagram-mod2.eps}}}
       \put(-14,2.5){\scalebox{0.5}{\includegraphics{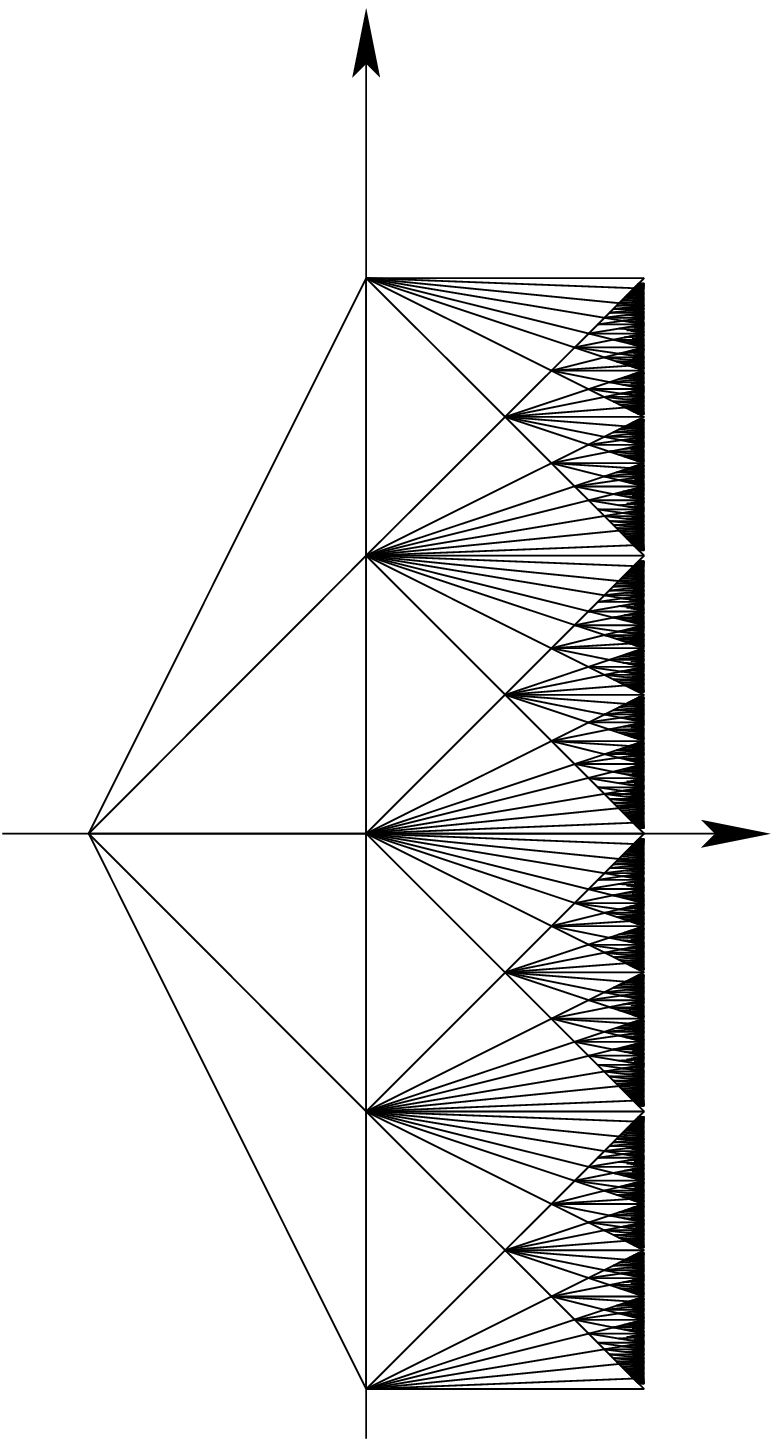}}}
       \put(58,175){\rotatebox{90}{\scalebox{1.0}{$\cdots$}}}
       \put(58,-10){\rotatebox{90}{\scalebox{1.0}{$\cdots$}}}
       \put(-27,108){$\angleb{1/0}$}
       \put(-13,104){\vector(1,-1){10}}
       \put(-14,83){$-1$}
%       \put(38,88){$0$}
       \put(82,82){$1$}
       \put(102,88){$u$}
       \put(45,205){$v$}
       \put(1,12){$\angleb{-2}$}
       \put(23,14){\vector(4,-1){15}}
       \put(26,3){$-2$}
       \put(-11,50){$\angleb{-1}$}
       \put(11,51){\vector(1,0){25}}
       \put(26,43){$-1$}
       \put(10,94){$\angleb{0}$}
       \put(23,95){\vector(4,-1){15}}
       \put(30,82){$O$}
       \put(-3,129){$\angleb{1}$}
       \put(11,131){\vector(1,0){25}}
       \put(33,134){$1$}
       \put(6,169){$\angleb{2}$}
       \put(20,171){\vector(1,0){15}}
       \put(33,173){$2$}
       \put(100,0){}
      \end{picture}
     \end{center}
     \caption{The diagram $\Diagram$}
     \label{Fig:Diagram}
%    \end{minipage}
  \end{figure}

The diagram $\Diagram$ is described as 
the 1-skeleton of a triangulation of a region on a plane 
as illustrated in Figure \ref{Fig:Diagram}.
%%% It is worth mentioning that
%The diagram $\Diagram$ is essentially the same as the so-called Farey graph.
%%% In addition, 
%%% since the term ``diagram'' sometimes means the diagram $\Diagram$
%%% and on another occasion a diagram of some knot or tangle,
%%% don't get confused between the two meanings.
%
Precisely a vertex of $\Diagram$ indicates either 
a point $(u,v)=((q-1)/q,p/q)$ denoted by $\angleb{ p/q }$, 
a point $(u,v)=(1,p/q)$ denoted by $\circleb{p/q}$
where $p/q$ is an irreducible fraction,
or a point $(u,v)=(-1,0)$ denoted by $\angleb{1/0}=\angleb{\infty}$.

  \begin{figure}[hbt]
%\hspace{-10mm}
%
%    \begin{minipage}{200pt}
%    \begin{minipage}{180pt}
    \begin{center}    
     \begin{picture}(100,85)
      \scalebox{0.75}{
       \put(0,-10){\scalebox{0.75}{\includegraphics{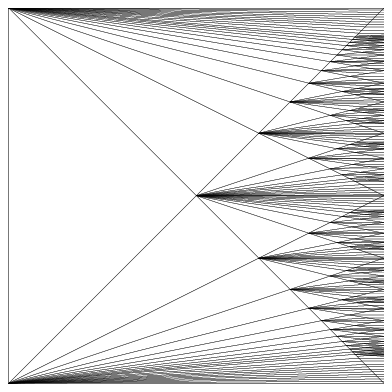}}}
      }
      \put(-9,-1){\vector(1,0){10}}
      \put(-23,-4){$\angleb{0}$}
      \put(121,-1){\vector(-1,0){10}}
      \put(124,-4){$\circleb{0}$}
      \put(-9,100){\vector(1,0){10}}
      \put(-23,97){$\angleb{1}$}
      \put(121,100){\vector(-1,0){10}}
      \put(124,97){$\circleb{1}$}
      \put(43,50){\vector(1,0){10}}
      \put(28,47){$\angleb{\frac{1}{2}}$}
      \put(121,50){\vector(-1,0){10}}
      \put(124,47){$\circleb{\frac{1}{2}}$}
      \put(61,33){\vector(1,0){10}}
      \put(46,30){$\angleb{\frac{1}{3}}$}
      \put(121,33){\vector(-1,0){10}}
      \put(124,30){$\circleb{\frac{1}{3}}$}
      \put(61,66){\vector(1,0){10}}
      \put(46,63){$\angleb{\frac{2}{3}}$}
      \put(121,66){\vector(-1,0){10}}
      \put(124,63){$\circleb{\frac{2}{3}}$}
     \end{picture}
    \end{center}
    \caption{
     A part of the diagram $\Diagram$ in 
     $[0,1]\times[0,1]$ 
    }
    \label{Fig:Diagram2}
%    \end{minipage}

%  \end{center}
% HAVE COPIED FROM LOWERBOUND PAPER
  \end{figure}

%
%	Edge
%

\subsubsection{Edges}

Two vertices $\angleb{p/q}$ and $\angleb{r/s}$ in $\Diagram$ 
are connected by an \textit{edge} if $|ps-qr|=1$. 
Such an edge is denoted by $\angleb{p/q}$\EDGE$\angleb{r/s}$.
An important class of the edges are the {\em vertical edges},
which connect the vertices 
$\angleb{z}$ and $\angleb{z+1}$ for arbitrary integer $z$.
%%% There is a vertex labeled by $\infty$ which corresponds to $1/0$.
Another important class of the edges are the {\em $\infty$-edges},
which connect the vertices 
$\angleb{\infty}$ and $\angleb{z}$ for integer $z$.
There is another kind of edge called a {\em horizontal edge},
which connects 
%%% $((q-1)/q,p/q)$ and $(1,p/q)$.
$\angleb{p/q}$ and $\circleb{p/q}$ for each $p/q$.

Let $e$ be an edge $\angleb{p/q}$\EDGE$\angleb{r/s}$
with $q\ge 1$, $s\ge 1$.
Let $k$ and $m$ be integers 
satisfying $m\ge 2$ and $1\le k \le m-1$.
Then,
let $\frac{k}{m}\angleb{p/q}+\frac{m-k}{m}\angleb{r/s}$
denote a point on $e$
with $uv$-coordinates $(u,v)=
\frac{kq}{kq+(m-k)s}(\frac{q-1}{q},\frac{p}{q})
+\frac{(m-k)s}{kq+(m-k)s}(\frac{s-1}{s},\frac{r}{s})$.
This is called a {\em rational point} of the edge $e$.
%
%%% Though the coordinates are 
%%% different from the usual internally dividing point
%%% $(u,v)=
%%% \frac{k}{m}(\frac{q-1}{q},\frac{p}{q})
%%% +\frac{m-k}{m}(\frac{s-1}{s},\frac{r}{s})$ in general,
%%% these concrete coordinates are not important 
%%% in the later argument.
%
For an $\infty$-edge $\angleb{1/0}$\EDGE$\angleb{z}$,
let
$\frac{k}{m}\angleb{1/0}+\frac{m-k}{m}\angleb{z}$
denote a point
with $(u,v)=\frac{k}{m}(-1,0)+\frac{m-k}{m}(0,z)$ on the edge,
which is a rational point of the $\infty$-edge.
A {\em partial edge} is an edge
$(\frac{k}{m}\angleb{p/q}+\frac{m-k}{m}\angleb{r/s})$\EDGE$\angleb{r/s}$.
%For a non-constant edgepath,
%the last edge can be a partial edge.
On the other hand,
a non-partial edge is called a {\em complete edge}.

Let $e$ be a horizontal edge $\angleb{p/q}$\EDGE$\circleb{p/q}$.
Then,
for integers $l\ge 1$ and $m\ge 1$
let $\frac{m}{m+l}\angleb{p/q}+\frac{l}{m+l}\circleb{p/q}$
denote a point
with $(u,v)=\frac{m}{m+l}(\frac{q-1}{q},\frac{p}{q})+\frac{l}{m+l}(1,\frac{p}{q})$ on the edge $e$.
This is a rational point of the horizontal edge.

Assume that
a non-$\infty$-edge is oriented from right to left. 
%an edgepath for $R$-tangle 
%is directed from $\angleb{R}$ to the other end.
%
Then  it is said to be {\em increasing} or {\em decreasing}
if the $v$-coordinate increases or decreases
as a point moves along the edge in that direction.
That is, an edge in the region $u>0$ is increasing or decreasing
if the $v$-coordinate increases or decreases 
as the $u$-coordinate decreases.

We assign $+1$ or $-1$ to an edge $e$ 
according to whether $e$ is increasing or decreasing respectively.
This is called the {\em sign} of the edge $e$,
and is denoted by $\sigma(e)$.

The \textit{length} of a complete edge is set to be $1$.
%%% For type I edgepaths, we have to deal with partial edges,
%%% which have length between $0$ and $1$. 
The length of a partial edge $(\frac{k}{m}\angleb{p/q}+\frac{m-k}{m}\angleb{r/s})$\DIREDGE$\angleb{r/s}$ is 
set to be $\frac{k}{m}$. 
%%%,
%%% but we omit the details here. 
Let $|e|$ denote the length of an edge $e$.

%
%	Edgepath
%

\subsubsection{Edgepaths}

Roughly, an edgepath is a connected piecewise-linear path in the diagram $\Diagram$. 
There are two kinds of edgepaths; 
non-constant edgepath and constant edgepath.

%Let $R$ be a fixed rational tangle.
%
The first one is a path on $\Diagram$ starting from a vertex $\angleb{R}$.
We will only consider such an edgepath with the following two conditions: 
\begin{enumerate}
\item
It must run from right to left in the weak sense, 
where vertical edges are allowed. 
\item
It must be {\em minimal}, that is,
it must not retrace back nor include two successive edges on a common triangle.
\end{enumerate}
We call an edgepath of this kind a {\em non-constant edgepath}.
A non-constant edgepath is expressed by a sequence of vertices
like
$\angleb{p_k/q_k}$\DIREDGE%
$\angleb{p_{k-1}/q_{k-1}}$\DIREDGE%
$\ldots$\DIREDGE%
$\angleb{p_2/q_2}$\DIREDGE%
$\angleb{p_1/q_1}$,
where $p_1/q_1=R$.
%%% Strictly speaking, a partial edge can be combined at the tail
%%% as mentioned later.
Note here that
the vertices are listed from right to left
so that the direction of the sequence coincides with 
the direction of an edgepath in the diagram $\Diagram$.
%
%%% and going right to left in the weak sense.
%%% An edgepath .
%
Also note that we allow that %For a non-constant edgepath,
the last edge can be a partial edge. 

An edgepath of the second kind
is a single point on a horizontal edge $\angleb{R}$\EDGE$\circleb{R}$. 
This is called a {\em constant edgepath}.
The single point of a constant edgepath 
must be a rational point of the horizontal edge
$\angleb{R}$\EDGE$\circleb{R}$.
%
%Though, constant edgepaths play not so important roles in this paper.

We will always use the symbol $\Edgepath$ for edgepaths.
In the following, we collect some terminologies about edgepaths.

%Let $\Edgepath$ be an edgepath.
Assume that the $u$-coordinate of the left endpoint of an edgepath is $u_0$.
Then the edgepath is said to be of {\em type I}, {\em type II} or {\em type III} 
according to which of $u_0>0$, $u_0=0$ or $u_0<0$ the coordinate satisfies.

We here introduce a new notion; a {\em basic edgepath}, which is not used in \cite{HO}. 
An edgepath is called a {\em basic edgepath} 
if the edgepath ends at the point with $u$-coordinate $0$ and includes no vertical edges.
We usually use symbol $\BasicEdgepath$ for basic edgepaths.

An edgepath is said to be {\em monotonically increasing} (respectively {\em monotonically decreasing}),
if all non-$\infty$-edges in the edgepath are increasing (resp. decreasing).
We ignore $\infty$-edges when we judge an edgepath is monotonic or not.
%
%
%%% For the vertex with label $p/q (q\ge 2)$,
For the vertex $\angleb{p/q}$ $(q\ge 2)$, 
there is one increasing leftward edge and one decreasing leftward edge. 
Hence, there exists just one monotonically increasing basic edgepath
and just one monotonically decreasing basic edgepath 
for a fixed fraction $p/q$ $(q\ge 2)$.

The \textit{length} of an edgepath is the sum of lengths of edges in the edgepath.
Naturally the length of a constant edgepath should be zero.
For an edgepath $\Edgepath$, $|\Edgepath|$ denotes the length of the edgepath.

%
%	Examples
%

\subsubsection{Examples}

\begin{figure}[hbt]
 \begin{picture}(300,390)
%1-3-3-1-2-b.eps PS 290x182 290x182+0+0 PseudoClass 256c 52kb 0.010u 0:01
  \put(0,200){
%   \put(0,0){\scalebox{1.3}{\includegraphics{edgepath-basic-1.eps}}}
   \put(0,0){\scalebox{1.3}{\includegraphics{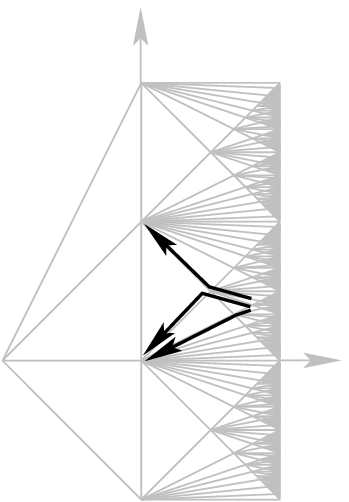}}} 
   % thick color c0c0c0
   \put(65,115){\vector(-1,-1){10}}
   \put(67,112){$\angleb{1}$}
   \put(95,90){\vector(-2,-1){15}}
   \put(97,87){$\angleb{1/2}$}
   \put(116,75){\vector(-1,0){20}}
   \put(118,72){$\angleb{2/5}$}
   \put(104,62){\vector(-2,1){15}}
   \put(106,59){$\angleb{1/3}$}
   \put(65,43){\vector(-1,1){10}}
   \put(67,40){$\angleb{0}$}
  }
  \put(50,192){(1)}
  \put(150,200){
%   \put(0,0){\scalebox{1.3}{\includegraphics{edgepath-basic-2.eps}}}
   \put(0,0){\scalebox{1.3}{\includegraphics{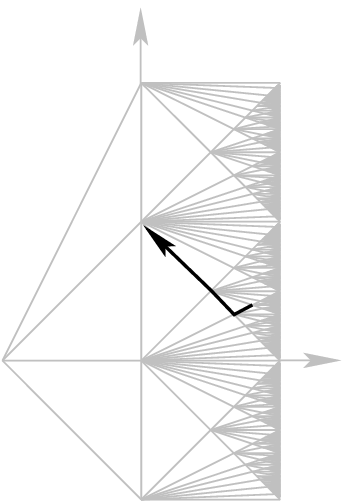}}}
   \put(65,115){\vector(-1,-1){10}}
   \put(67,112){$\angleb{1}$}
   \put(95,90){\vector(-2,-1){15}}
   \put(97,87){$\angleb{1/2}$}
   \put(116,75){\vector(-1,0){20}}
   \put(118,72){$\angleb{2/5}$}
   \put(104,62){\vector(-2,1){15}}
   \put(106,59){$\angleb{1/3}$}
  }
  \put(200,192){(2)}
  \put(0,0){
%   \put(0,0){\scalebox{1.3}{\includegraphics{edgepath-type13-mod.eps}}}
   \put(0,0){\scalebox{1.3}{\includegraphics{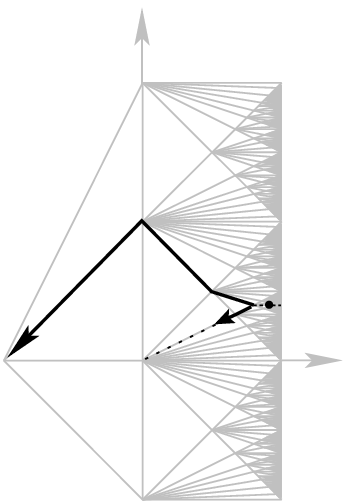}}}
   \put(3,70){\vector(0,-1){15}}
   \put(-9,76){$\angleb{1/0}$}
   \put(65,115){\vector(-1,-1){10}}
   \put(67,112){$\angleb{1}$}
   \put(80,95){\vector(0,-1){15}}
   \put(70,99){$\angleb{1/2}$}
   \put(116,95){\vector(-1,-1){20}}
   \put(118,92){$\angleb{2/5}$}
%   \put(121,87){\vector(-2,-1){23}}
%   \put(123,84){$\angleb{2/5}$}
   \put(116,75){\vector(-1,0){10}}
   \put(118,72){$\circleb{2/5}$}
   \put(104,62){\vector(-2,1){15}}
   \put(106,59){$\angleb{1/3}$}
   \put(65,43){\vector(-1,1){10}}
   \put(67,40){$\angleb{0}$}
  }
  \put(50,-8){(3)}
  \put(150,0){
%   \put(0,0){\scalebox{1.3}{\includegraphics{edgepath-type2.eps}}}
   \put(0,0){\scalebox{1.3}{\includegraphics{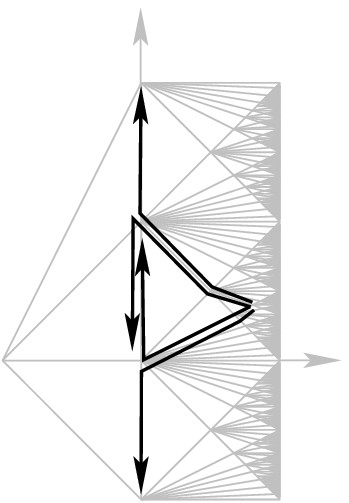}}}
   \put(65,145){\vector(-1,1){10}}
   \put(67,142){$\angleb{2}$}
   \put(65,115){\vector(-1,-1){10}}
   \put(67,112){$\angleb{1}$}
   \put(95,90){\vector(-2,-1){15}}
   \put(97,87){$\angleb{1/2}$}
   \put(116,75){\vector(-1,0){20}}
   \put(118,72){$\angleb{2/5}$}
   \put(104,62){\vector(-2,1){15}}
   \put(106,59){$\angleb{1/3}$}
   \put(65,43){\vector(-1,1){10}}
   \put(67,40){$\angleb{0}$}
   \put(65,13){\vector(-1,-1){10}}
   \put(67,10){$\angleb{-1}$}
  }
  \put(200,-8){(4)}
 \end{picture}
 \caption{Examples of edgepaths for $2/5$-tangle}
 \label{Fig:Edgepaths}
\end{figure}

Some examples of edgepaths for $2/5$-tangle are given in Figure \ref{Fig:Edgepaths}.
In (1), three basic edgepaths 
$\angleb{1}$\DIREDGE$\angleb{1/2}$\DIREDGE$\angleb{2/5}$,
$\angleb{0}$\DIREDGE$\angleb{1/2}$\DIREDGE$\angleb{2/5}$ and
$\angleb{0}$\DIREDGE$\angleb{1/3}$\DIREDGE$\angleb{2/5}$
are drawn.
The first and third one are monotonically increasing and decreasing, while the second one is neither. 
In (2),
the edgepath
$\angleb{1}$\DIREDGE$\angleb{1/2}$\DIREDGE$\angleb{1/3}$\DIREDGE$\angleb{2/5}$
is not minimal and is not allowed,
since $\angleb{1/2}$\DIREDGE$\angleb{1/3}$\DIREDGE$\angleb{2/5}$
consists of two edges of a common triangle.
In (3),
a type I edgepath
($\frac{1}{2}\angleb{0}$+$\frac{1}{2}\angleb{1/3}$)\DIREDGE$\angleb{1/3}$\DIREDGE$\angleb{2/5}$
and
a type III edgepath
$\angleb{1/0}$\DIREDGE$\angleb{1}$\DIREDGE$\angleb{1/2}$\DIREDGE$\angleb{2/5}$
are drawn.
A thick dot on the horizontal edge
$\angleb{2/5}$\EDGE$\circleb{2/5}$
in this figure
is a constant edgepath,
which is a rational point on the horizontal edge.
%%% $\angleb{2/5}$\EDGE$\circleb{2/5}$.
In (4),
type II edgepaths
$\angleb{2}$\DIREDGE$\angleb{1}$\DIREDGE$\angleb{1/2}$\DIREDGE$\angleb{2/5}$,
$\angleb{0}$\DIREDGE$\angleb{1}$\DIREDGE$\angleb{1/2}$\DIREDGE$\angleb{2/5}$,
$\angleb{1}$\DIREDGE$\angleb{0}$\DIREDGE$\angleb{1/3}$\DIREDGE$\angleb{2/5}$ and
$\angleb{-1}$\DIREDGE$\angleb{0}$\DIREDGE$\angleb{1/3}$\DIREDGE$\angleb{2/5}$
are depicted.
Note that the second one is also an example of a non-minimal edgepath,
since the edgepath is not minimal at the vertex $\angleb{1}$.

%

%
%	Edgepath systems
%

\subsubsection{Edgepath systems}

By an {\em edgepath system}, we mean a collection of $\NumTangles$ edgepaths.
Recall that $\NumTangles$ denotes the number of tangles of a Montesinos knot. 
We will use symbol $\EdgepathSystem$ for edgepath systems. 

%For the consistency on gluing the subsurfaces corresponding to edgepaths in an edgepath system,
%edgepaths must satisfy the following conditions:
We call the following conditions {\em gluing consistency}: 

\smallskip

The $u$-coordinates of the left endpoints of the edgepaths must coincide with each other,
and the sum of $v$-coordinates of the left endpoints must be $0$.

\smallskip

As we will see, this condition must be satisfied by an edgepath system 
which corresponds to a properly embedded surface in the exterior of a Montesinos knot. 

%Note that the left endpoint of an edgepath 
%expresses the boundary of a subsurface corresponding to the edgepath,
%which is a curve system.

Edgepaths in an edgepath system satisfying the gluing consistency 
have the common $u$-coordinate of the left endpoints.
Therefore, similarly to edgepaths,
edgepath systems satisfying the gluing consistency 
are also classified by the common $u$-coordinates.
%%% of the left endpoints.
That is,
such an edgepath system is said to be of {\em type I}, {\em type II} or {\em type III}
according to which of $u>0$, $u=0$ or $u<0$ the common $u$-coordinate satisfies.

In general, we will only consider edgepath systems with gluing consistency, 
%In the algorithm,
%an essential surface corresponds to
%an edgepath system satisfying gluing consistency.
%Though,
but we sometimes make an edgepath system 
which may not satisfy the gluing consistency.
We call the edgepath system a {\em formal edgepath system}.
%
%Of course,
%a formal edgepath system may not correspond to an essential surface.
%However, it does not matter,
%since ``twist'' mentioned later can be calculated also for formal edgepath systems,
%and these edgepath systems are used only to give some value of ``twist''.

For example, 
we will consider an edgepath system all of whose edgepaths are basic edgepaths, 
which we call a {\em basic edgepath system}. 
Such a basic edgepath system is formal. 
We usually use symbol $\BasicEdgepathSystem$ for basic edgepath systems.

\subsubsection{Diagram and curve systems}

We here observe a relationship between 
vertices and edges in $\Diagram$ and curve systems on a fourth-punctured sphere. 
This enables us to describe a connection of edgepath systems and 
embedded surfaces in Montesinos knot exteriors. 
Such a connection will be explained in the next subsubsection. 

\begin{figure}[hbt]
 \begin{picture}(140,75)
  \put(0,12){\scalebox{0.6}{\includegraphics{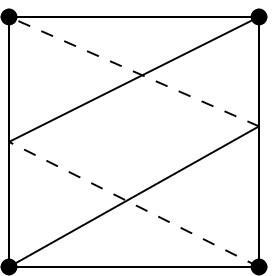}}}
  \put(80,0){\scalebox{0.6}{\includegraphics{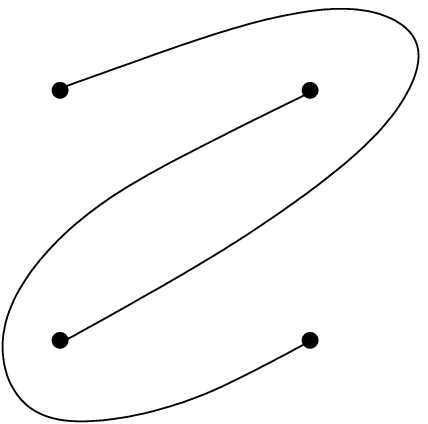}}}
 \end{picture}
 \caption{$1/2$-arcs on a flattened sphere and on a plane.}
 \label{Fig:1over2-Arcs}
\end{figure}

\begin{figure}[hbt]
 \begin{picture}(220,60)
%1-3-3-1-2-b.eps PS 290x182 290x182+0+0 PseudoClass 256c 52kb 0.010u 0:01
  \put(0,0){\scalebox{0.6}{\includegraphics{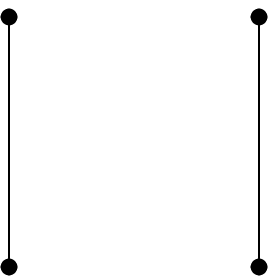}}}
  \put(80,0){\scalebox{0.6}{\includegraphics{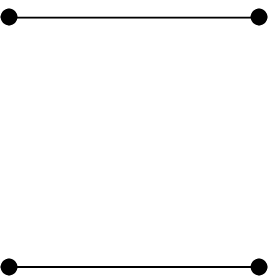}}}
  \put(160,0){\scalebox{0.6}{\includegraphics{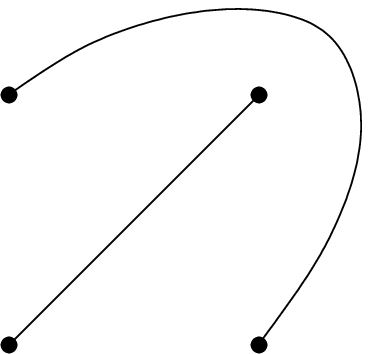}}}
 \end{picture}
 \caption{$1/0$-arcs, $0$-arcs, $1$-arcs}
 \label{Fig:SomePairsOfArcs}
\end{figure}

As mentioned before,
segments of slope $p/q$ on a pillowcase form two arcs.
We call the arcs {\em (a pair of) $p/q$-arcs}.
For instance,
Figure \ref{Fig:1over2-Arcs} (left) depicts $1/2$-arcs.
For ease of drawing figures,
we regard $\mathbb{R}^2\cup\{\infty\}$ as a sphere
and draw the four punctures and $p/q$-arcs on the plane.
See $1/2$-arcs in Figure \ref{Fig:1over2-Arcs} (right).
Some other examples of $p/q$-arcs are given in Figure \ref{Fig:SomePairsOfArcs}.
We define that the set of $m$ pairs of $p/q$-arcs correspond to the vertex $\angleb{p/q}$ on $\Diagram$. 
Please note that such multiple pairs 
can be drawn on the fourth-punctured sphere disjointly.

Next, consider the integers $p$, $q$, $r$ and $s$ satisfying $|ps-qr|=1$. 
Then a combination of $m$ pairs of $p/q$-arcs and $m-k$ pairs of $r/s$-arcs 
is set to correspond to a rational point $\frac{k}{m}\angleb{p/q}+\frac{m-k}{m}\angleb{r/s}$.
Also note that such a pair of $p/q$-arcs and a pair of $r/s$-arcs 
can be disjointly drawn on a fourth-punctured sphere.
A disjoint union of such pairs appears in Figure \ref{Fig:SubsurfaceMultiSheets}.

%
%%% Strictly speaking,
%%% we can construct a disjoint union 
%%% of pairs of $p_1/q_1$-arcs, pairs of $p_2/q_2$-arcs and pairs of $p_3/q_3-arcs
%%% for fractions $p_1/q_1$, $p_2/q_2$ and $p_3/q_3$
%%% for which $\angleb{p_1/q_1}$, $\angleb{p_2/q_2}$ and $\angleb{p_3/q_3}$
%%% form a triangle in $\diagram$.
%%% But this kind of pattern does not appear in the Hatcher-Oertel algorithm.
%

Besides,
a non-null-homotopic loop on a fourth-punctured sphere disjoint from $p/q$-arcs
is topologically unique and is called a {\em $p/q$-circle}.
We can think a disjoint union of $m$ pairs of $p/q$-arcs and $l$ copies of $p/q$-circles with $m\ge 1$, $l\ge 1$,
as seen in Figure \ref{Fig:SubsurfaceCap}.
This corresponds to a rational point $\frac{m}{m+l}\angleb{p/q}+\frac{l}{m+l}\circleb{p/q}$ of a horizontal edge.

\subsubsection{Edgepath systems and corresponding surfaces}

Consider the same setting as in the explanation 
of the outline of Hatcher-Oertel's machinery at the beginning of this subsection. 
Then a subsurface $F_i$ lies in one of the $N$ $3$-balls.
Now,
we assume that the $3$-ball is actually a unit ball
with the center $O$.
We think about the intersection between $F_i$
and a sphere with the center $O$ and radius $\Radius>0$.
As we vary $\Radius$ within $0<\Radius\le 1$,
the intersection changes.
The subsurface $F_i$ is expressed by the sequence of
such intersections between $F_i$ and concentric spheres.

%In the algorithm,
By virtue of Proposition 1.1 in \cite{HO}, 
%a subsurface is isotoped in advance
%into some normalized form.
%In particular,
the subsurface is deformed into a good position
in Morse theoretical sense.
This enables us to regard
the intersection as a curve system,
on any non-critical level $\Radius$. 
The topological type of the intersection changes
only at critical levels, 
and so, 
a subsurface can be represented by a sequence of curve systems
and hence by an edgepath.

In the following, we observe parts of the surfaces 
represented by each edge in the edgepath system. 

%\noindent {(Saddle subsurface)}

\begin{figure}[hbt]
 \begin{picture}(380,90)
%1-3-3-1-2-b.eps PS 290x182 290x182+0+0 PseudoClass 256c 52kb 0.010u 0:01
  \put(0,30){\scalebox{0.6}{\includegraphics{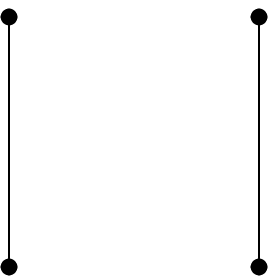}}}
  \put(-7,0){$3/4< \Radius\le 1$}
  \put(80,30){\scalebox{0.6}{\includegraphics{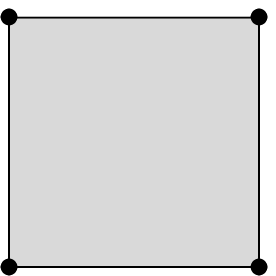}}}
  \put(83,0){$\Radius=3/4$}
  \put(160,30){\scalebox{0.6}{\includegraphics{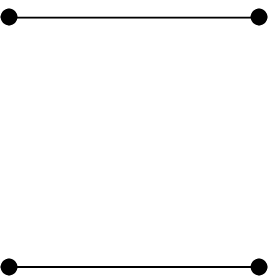}}}
  \put(150,0){$1/2< \Radius<3/4$}
  \put(240,17){\scalebox{0.6}{\includegraphics{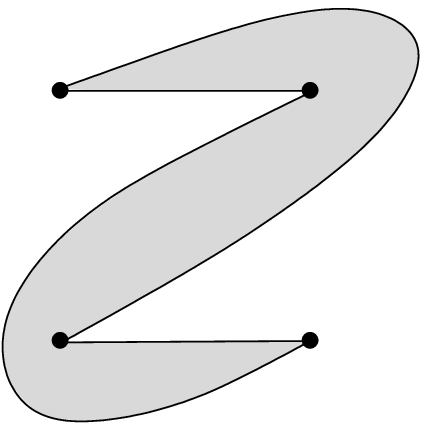}}}
  \put(243,0){$\Radius=1/2$}
  \put(320,17){\scalebox{0.6}{\includegraphics{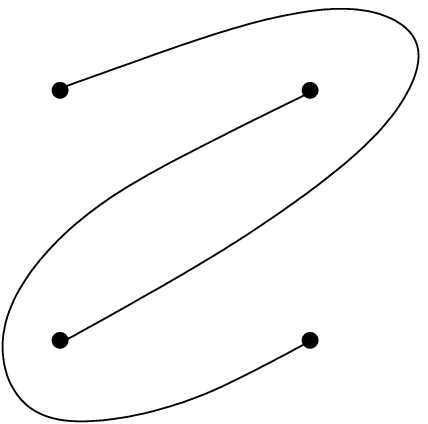}}}
  \put(310,0){$1/4\le \Radius<1/2$}
 \end{picture}
 \caption{The intersection between a concentric sphere of radius $\Radius$
          and a subsurface
          expressed by an edgepath
          $\angleb{1/0}$\DIREDGE$\angleb{0}$\DIREDGE$\angleb{1/2}$.}
 \label{Fig:SubsurfaceSaddle}
\end{figure}

A non-constant edgepath corresponds to a {\em saddle subsurface}.
An example of a saddle subsurface
is expressed as a sequence of intersections with concentric spheres
as illustrated in Figure \ref{Fig:SubsurfaceSaddle}.
In this figure,
the intersection is 
empty for $\Radius<1/4$,
$1/2$-arcs for $1/4\le \Radius<1/2$,
a disk bounded by $1/2$-arcs and $0$-arcs for $\Radius=1/2$,
$0$-arcs for $1/2<\Radius<3/4$,
a disk bounded by $0$-arcs and $1/0$-arcs for $\Radius=3/4$,
and 
$1/0$-arcs for $3/4<\Radius\le 1$.
%%% (Of course, the exact value of threshold of $\Radius$ like $1/4$, $1/2$ and $3/4$
%%% is not important topologically.)
%
%
This is a subsurface corresponding
to an edgepath $\angleb{1/0}$\DIREDGE$\angleb{0}$\DIREDGE$\angleb{1/2}$.
A disk bounded by two pairs of arcs is called a {\em saddle}.
A pair of $p/q$-arcs and a pair of $r/s$-arcs
can be disjointly drawn and give a saddle only if $|ps-qr|=1$.
This is the reason why we think about an edgepath 
which consists of edges of the form $\angleb{p/q}$\EDGE$\angleb{r/s}$
satisfying $|ps-qr|=1$.
%
%%% Strictly speaking,
%%% there is an ambiguity on making a saddle.
%%% That is,
%%% the two pairs of arcs on a sphere form a loop,
%%% there exist two disks on a sphere bounded by the loop,
%%% and hence there are two choices of making a saddle.
%%% Thus,
%%% there remains some ambiguity 
%%% on making a subsurface from an edgepath.
%%% Though,
%%% this ambiguity is not important in this paper
%%% and we can ignore it,
%%% since all subsurfaces sharing the common edgepath 
%%% contribute the same amount to the boundary slope.

Note that,
four punctures on each concentric sphere with radius $1/4<\Radius\le 1$
and $1/2$-arcs on the sphere with radius $\Radius=1/4$
form a $1/2$-tangle in the $3$-ball.
Thus, the boundary of the subsurface
lies on the boundary of the exterior of $1/2$-tangle.

%\noindent (The Number of sheets)

The example of a subsurface above is ``one-sheeted''.
In general, a subsurface can have ``multiple sheets''.
If a small meridian circle of a Montesinos knot $K$ meets 
an essential surface in $m$ points, then the surface is said to be {\em $m$-sheeted}, 
or the {\em number of sheets} of the surface is $m$ at the point on the knot. 
After isotoping the surface with keeping it properly embedded,
%%%  (or $3$-ball including the surface)
the number of sheets is constant everywhere on the knot.
Then, the number of sheets is defined for a surface.
For an $m$-sheeted subsurface $F$,
the intersection between $F$ and a concentric sphere
is in general a curve system including exactly $m$ pairs of arcs.
%%% where the union can include
%%% some null-homotopic circles if $F$ is a cap subsurface.

%\noindent {(Saddle subsurface with multiple sheets)}

\begin{figure}[hbt]
 \begin{picture}(250,75)
  \put(0,15){\scalebox{0.6}{\includegraphics{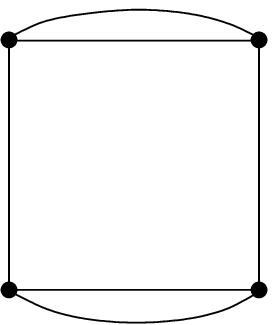}}}
  \put(-7,0){$2/3< \Radius\le 1$}
  \put(100,15){\scalebox{0.6}{\includegraphics{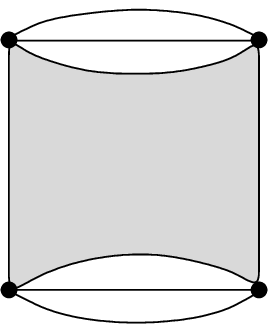}}}
  \put(101,0){$\Radius=2/3$}
  \put(200,15){\scalebox{0.6}{\includegraphics{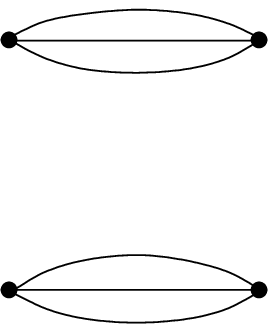}}}
  \put(190,0){$1/3\le \Radius<2/3$}
 \end{picture}
 \caption{The intersection between a concentric sphere of radius $\Radius$
          and a subsurface
          expressed by an edgepath
          $(\frac{1}{3}\angleb{1/0}+\frac{2}{3}\angleb{0})$\DIREDGE$\angleb{0}$.}
 \label{Fig:SubsurfaceMultiSheets}
\end{figure}

%
Figure \ref{Fig:SubsurfaceMultiSheets} illustrates an example of a saddle subsurface with multiple sheets.
In this figure, there are $3$ sheets, and only one saddle appears.
The left figure is the union of one pair of $1/0$-arcs and two pairs of $0$-arcs,
which is denoted by $\frac{1}{3}\angleb{1/0}+\frac{2}{3}\angleb{0}$.
Then, this subsurface corresponds to a partial edge
$(\frac{1}{3}\angleb{1/0}+\frac{2}{3}\angleb{0})$\DIREDGE$\angleb{0}$.
An edgepath including a partial edge must be multiple-sheeted.

Note that, since $\angleb{p/q}$ corresponds also to
a disjoint union of $m$ pairs of $p/q$-arcs for arbitrary positive integer $m$,
an $m$-sheeted parallel subsurface also corresponds to an edgepath without a partial edge.

%\noindent {(Cap subsurface)}

\begin{figure}[hbt]
 \begin{picture}(280,80)
  \put(0,15){\scalebox{0.6}{\includegraphics{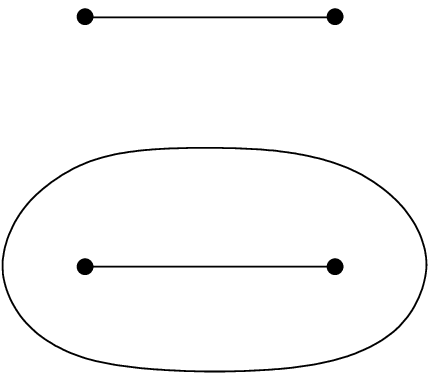}}}
  \put(5,0){$2/3\le \Radius\le 1$}
  \put(100,15){\scalebox{0.6}{\includegraphics{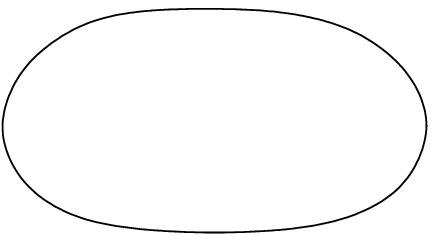}}}
  \put(101,0){$1/3<\Radius<2/3$}
  \put(200,15){\scalebox{0.6}{\includegraphics{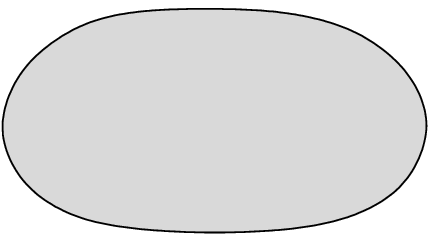}}}
  \put(213,0){$\Radius=1/3$}
 \end{picture}
 \caption{The intersections between a concentric sphere of radius $\Radius$
          and a subsurface
          expressed by an edgepath
          $(\frac{1}{2}\angleb{0}+\frac{1}{2}\circleb{0})$.}
 \label{Fig:SubsurfaceCap}
\end{figure}

%
A subsurface corresponding to a constant edgepath is illustrated in Figure \ref{Fig:SubsurfaceCap}.
The intersection is empty for $\Radius<1/3$,
a disk bounded by a $0$-circle for $\Radius=1/3$,
a $0$-circle for $1/3<\Radius<2/3$,
and
the union of a $0$-circle and $0$-arcs for $2/3\le \Radius\le 1$.
Such a subsurface is called a {\em cap subsurface}.
Note that the boundary of this subsurface 
lies on the boundary of the exterior of $0$-tangle in the $3$-ball.
%%% similarly to the previous example.
%
The intersection between the subsurface and a concentric sphere of radius $2/3\le \Radius\le 1$
is expressed by a rational point $\frac{1}{2}\angleb{0}+\frac{1}{2}\circleb{0}$.
And the subsurface is also denoted by
$\frac{1}{2}\angleb{0}+\frac{1}{2}\circleb{0}$, which means a constant edgepath.
%%%,
%%%which is regarded as a one-point edgepath called a constant edgepath.
%%% This is because
%%% we can recover the cap subsurface
%%% from $\frac{1}{2}\angleb{0}+\frac{1}{2}\circleb{0}$.
%
%
%
%

%\noindent {(Boundaries of subsurfaces and surface)}

For a saddle subsurface or a cap subsurface,
the boundary of the subsurface
lies on the boundary of the exterior of the corresponding rational tangle.
Therefore,
after the subsurfaces are glued each other,
the boundary of the obtained surface
lies on the boundary of the exterior of a Montesinos knot.
In short, the surface is confirmed to be properly embedded.

\subsubsection{Twists and boundary slopes}

\begin{figure}[hbt]
 \begin{picture}(200,86)
  \put(0,22){\scalebox{0.6}{\includegraphics{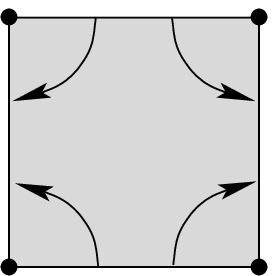}}}
  %\put(-7,0){$2/3< \Radius\le 1$}
  \put(100,0){\scalebox{0.6}{\includegraphics{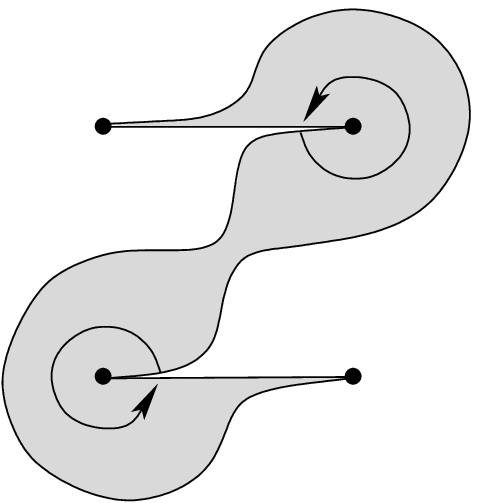}}}
  %\put(101,0){$\Radius=2/3$}
 \end{picture}
 \caption{Saddles for $\angleb{1/0}$\DIREDGE$\angleb{0}$
          and for $\angleb{0}$\DIREDGE$\angleb{1/2}$.}
 \label{Fig:Twist}
\end{figure}

In order to calculate the boundary slope of an essential surface,
we calculate the {\em total number of twists}, or in short {\em twist}.
As well as the boundary slope,
the twist also measures how many times
the boundary of a surface runs around a knot in the counterclockwise direction of the meridian,
while the boundary runs once along the knot.
%%% in anti-clockwise direction.
The twist of a surface is obtained
by first calculating the twist for each subsurface
and then summing up the twists.
The twist of a subsurface
is calculated as the sum of twists for each saddle in the subsurface.
For a saddle corresponding to an $\infty$-edge like $\angleb{1/0}$\DIREDGE$\angleb{z}$,
the twists appear on a saddle at four punctures
and cancel out each other as in Figure \ref{Fig:Twist} (left).
Besides,
for a saddle corresponding to a non-$\infty$-edge $\angleb{p/q}$\DIREDGE$\angleb{r/s}$,
the twist is $2$ or $-2$
as seen in Figure \ref{Fig:Twist} (right).
The sign of the twist is determined by whether the edge is decreasing or increasing.
Next,
for saddles of a subsurfaces with multiple sheets
as in Figure \ref{Fig:SubsurfaceMultiSheets},
a single saddle of an $m$-sheeted subsurface
corresponding to
a non-$\infty$-edge 
$(\frac{k}{m}\angleb{p/q}+\frac{m-k}{m}\angleb{r/s})$\DIREDGE$(\frac{k+1}{m}\angleb{p/q}+\frac{m-k-1}{m}\angleb{r/s})$
has the twist $\pm 2/m$,
and the twist for $(\frac{k}{m}\angleb{p/q}+\frac{m-k}{m}\angleb{r/s})$\DIREDGE$\angleb{r/s}$
is $\pm 2k/m$.
For a cap subsurface as in Figure \ref{Fig:SubsurfaceCap},
since no saddles appear,
the twist
%%% for the cap subsurface
is $0$.

%%%By putting them all together,
Hence,
the procedure to calculate the twist is summarized as follows.
For each non-$\infty$-edge $e$ in an edgepath,
we assign $-2 \sigma(e) |e|$ to the edge $e$.
Then the twist of an edgepath is the sum of such assigned values 
for edges in the edgepath,
where a constant edgepath has twist $0$.
The twist of an edgepath system is the sum of twists 
of all edgepaths in the edgepath system.
The twist is a substitution of the boundary slope
and fits well with the algorithm.
We use the symbol $\Twist$ for the twist,
and
%%% $\Twist()$ 
$\Twist(e)$, 
$\Twist(\Edgepath)$, 
$\Twist(\EdgepathSystem)$ and $\Twist(F)$
denote the twist of 
an edge $e$,
an edgepath $\Edgepath$,
an edgepath system $\EdgepathSystem$ and a surface $F$.
Since the boundary slope is defined so that
the boundary slope of a Seifert surface $F_\Seifert$ is $0$,
the boundary slope of an essential surface $F$
is calculated as the difference of
the twist of the surface
and that of a Seifert surface,
that is,
$\Twist(F)-\Twist(F_\Seifert)$.

\subsubsection{Remarks}

We here collect some remarks 
about the previous subsubsections. 

\begin{itemize}
\item
Since each rational tangle $R_i$ in a Montesinos knot in this paper
is non-integral,
the starting point of an edgepath
must be a vertex $\angleb{p/q}$ or a rational point on a horizontal edge $\angleb{p/q}$\EDGE$\circleb{p/q}$ for a non-integral fraction $p/q$.
On the other hand,
Figures \ref{Fig:SubsurfaceMultiSheets} and \ref{Fig:SubsurfaceCap}
correspond to edgepaths
with the starting point on $\angleb{0}$\EDGE$\circleb{0}$.
Hence, 
these figures
are not appropriate in some sense.
These figures are used just for ease of drawing and explanation.

\item
%\noindent (Figures of subsurfaces in another style)
Here we illustrates three-dimensional figures of subsurfaces corresponding to 
a non-constant edgepath and a constant edgepath (See Figure \ref{Fig:PartialSurfaces}).
Note that 
these subsurfaces lie in $3$-balls
and
top boundaries of these subsurfaces lie on the boundary spheres of the $3$-balls,
which are not drawn in this figure.
Besides,
the ``cap'' of the cap subsurface in Figure \ref{Fig:PartialSurfaces}(b)
is somewhat deformed
from the subsurface illustrated in Figure \ref{Fig:SubsurfaceCap}.
\begin{figure}[hbt]
 \begin{picture}(330,120)
 \put(40,20){\scalebox{1.0}{\includegraphics{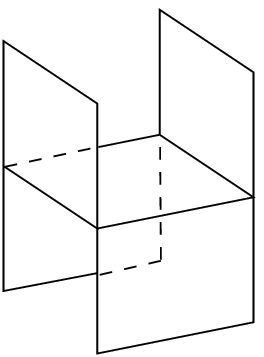}}} % 75x101
 \put(0,0){(a) A saddle subsurface}
 \put(18,-10){corresponding to $\angleb{\infty}$\DIREDGE$\angleb{0}$}
 \put(220,20){\scalebox{0.8}{\includegraphics{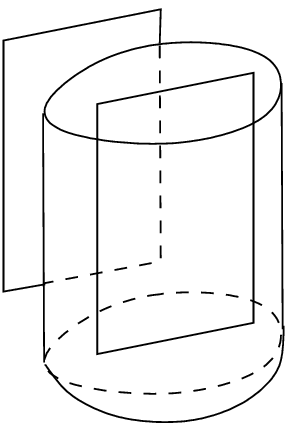}}} % 83x122
 \put(180,0){(b) A cap subsurface}
 \put(200,-10){corresponding to $\frac{1}{2}\angleb{0}+\frac{1}{2}\circleb{0}$} 
 \end{picture}
 \caption{Examples of subsurfaces}
 \label{Fig:PartialSurfaces}
\end{figure}

\item
Strictly speaking,
there are some irregular essential surfaces:
a Conway sphere and an augmented surface.
%%% The former does not correspond to an edgepath system,
%%% while the latter corresponds to an irregular edgepath system.
Although, we can ignore these surfaces
since their boundary slopes are $\infty$ or the same as that of some regular essential surface.
See \cite{HO} for detail.

\end{itemize}

%
% Subsection
%

\subsection{Crossing number of Montesinos knot}

In this subsection, we will explain the result given in \cite{LT} 
on the calculation of the crossing number for a Montesinos knot.

%
% Block
%

\subsubsection{Alternating tangle diagram}
%\ \newline\noindent {\bf Alternating tangle diagram}

Assume that $R$ is a non-integral fraction.
There is the unique continued fraction of $R$ as
\begin{eqnarray*}
 R=z_1+\frac{1}{z_2+\frac{1}{\ddots \frac{1}{z_{k-1}+\frac{1}{z_k}}}}
\end{eqnarray*}
with $k\ge 2$
where the integers $z_1$, $z_2$, $\ldots$, $z_k$ satisfy
\begin{eqnarray}
 \left\{
  \begin{array}{l}
   \textrm{ $z_1\ge 0$, $z_j\ge 1$ \,($1<j<k$), $z_k\ge 2$ \ \ if $R>0$},
   \\
   \\
   \textrm{ $z_1\le 0$, $z_j\le -1$ \,($1<j<k$), $z_k\le -2$ \ \ if $R<0$}.
  \end{array}
 \right.
\label{Eq:ConditionOfUniqueContinuedFraction}
\end{eqnarray}
Let $(z_1, z_2, \ldots, z_k)$ denote this continued fraction.
For instance, a continued fraction expansion of $47/36$ is
\begin{eqnarray*}
\frac{47}{36}&=&1+\frac{11}{36} 
=1+\frac{1}{\frac{36}{11}}=1+\frac{1}{3+\frac{3}{11}} 
=1+\frac{1}{3+\frac{1}{\frac{11}{3}}}
%% \\
%% &=&1+\frac{1}{3+\frac{1}{3+\frac{2}{3}}}
%% =1+\frac{1}{3+\frac{1}{3+\frac{1}{\frac{3}{2}}}}
=
\ldots
\\
&=&1+\frac{1}{3+\frac{1}{3+\frac{1}{1+\frac{1}{2}}}}
%% \\
=(1,3,3,1,2).
\end{eqnarray*}

\begin{figure}[hbt]
 \begin{picture}(145,91)
%1-3-3-1-2-b.eps PS 290x182 290x182+0+0 PseudoClass 256c 52kb 0.010u 0:01
  \put(0,0){\scalebox{0.5}{\includegraphics{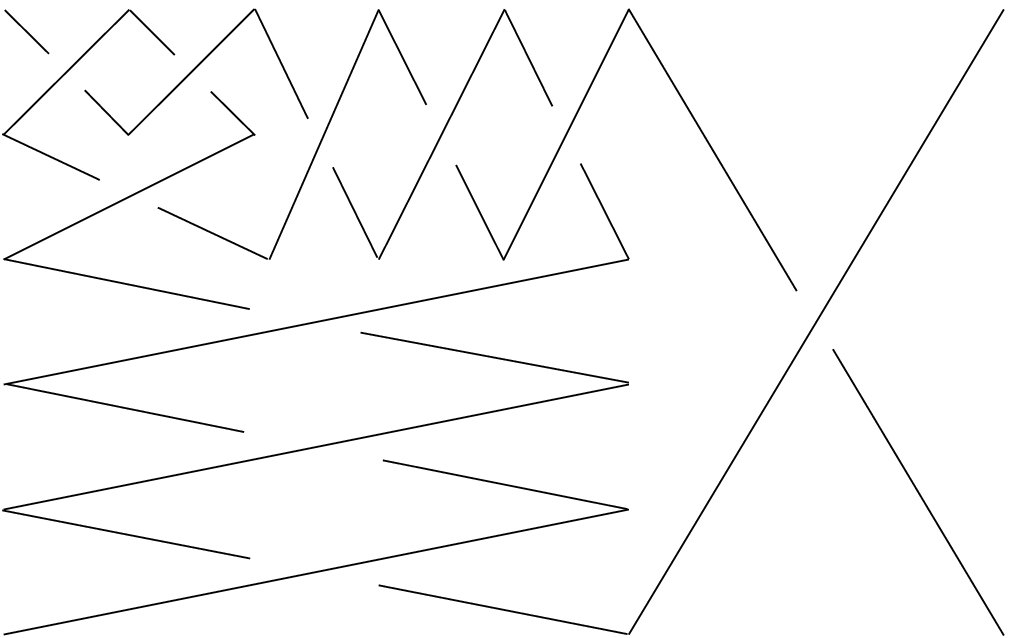}}}
 \end{picture}
 \caption{The alternating tangle diagram of the tangle $(1,3,3,1,2)=47/36$}
 \label{Fig:AlternatingDiagram47over36}
\end{figure}

This continued fraction gives a diagram of a rational tangle $R$,
which is roughly a combination of vertical and horizontal half-twists.
See Figure \ref{Fig:AlternatingDiagram47over36} for example.
If we connect two pairs of adjacent ends of the diagram,
because of the condition on $z_j$'s,
we obtain an alternating diagram of a two-bridge link.
Thus, we call a diagram like this of a rational tangle
an \textit{alternating tangle diagram} of a rational tangle $R$,
which is denoted by $D(R)$.
Note that the number of crossings of the diagram is
$\sum_{l=1}^{k} |z_l|$.
%
%
% Figure
%

%
% Block
%

\subsubsection{Reduced Montesinos diagram}
%\ \newline\noindent {\bf Reduced Montesinos diagram}

In \cite{LT}, Lickorish and Thistlethwaite determine the crossing number of the Montesinos knots.
Though, we must perform a technical operation which restricts the form of a diagram of a Montesinos knot.

Let $K$ be a Montesinos knot $K=K(R_1,R_2,\ldots,R_\NumTangles)$ with $\NumTangles \ge 3$, where each $R_j$ is a non-integral fraction.
For $(R_1,R_2,\ldots,R_\NumTangles)$,
let $D(R_1,R_2,\ldots,R_\NumTangles)$
denote the diagram of $K$ obtained by combining alternating tangle diagrams
$D(R_1)$, $D(R_2)$, $\ldots$, $D(R_\NumTangles)$ in a natural way.

Now,
for any tuple of integers $(i_1$, $i_2$, $\ldots$, $i_\NumTangles)$ 
with $i_1+i_2+\ldots+i_\NumTangles=0$,
the Montesinos knot $K^\prime=K(R_1+i_1,R_2+i_2,\ldots,R_\NumTangles+i_\NumTangles)$ is isotopic to $K$,
since addition of such integers corresponds to rotations of the tangles.
Thus,
for instance, if $\{R_j\}$ includes both positive and negative fractions and $R_a>1$ for some $a$,
by subtracting $1$ from $R_a$ and adding $1$ to some negative $R_b$,
we can reduce the integral part $[R_a]$ of $R_a$.
Through repeating such operations,
we finally have a tuple $(R^\prime_1, R^\prime_2, \ldots, R^\prime_\NumTangles)$
satisfying \vspace{1mm}\\
\hspace*{1cm} (A) all $R^\prime_j$'s have the same sign \vspace{1mm}\\
or \vspace{1mm} \\
\hspace*{1cm} (B) $\{R^\prime_j\}$ includes both positive and negative fractions, and $|R^\prime_j|<1$ holds for all $j$. \vspace{1mm} \\

For $(R_1,R_2,\ldots,R_\NumTangles)$ satisfying (A) or (B),
we think about the diagram $D=D(R_1,R_2,\ldots,R_\NumTangles)$.
If $\{ R_j \}$ have the common sign,
then the diagram satisfies condition (i) in \cite{LT},
%%% the corresponding Montesinos knot satisfies
that is, the diagram is alternating.
If $|R_j|<1$ holds for each $j$,
then each alternating tangle diagram $D(R_j)$ in the diagram $D$ satisfies
a technical condition in the condition (ii).
Therefore, $D$ satisfies the condition (i) or (ii) in \cite{LT}
to be a \textit{reduced Montesinos diagram}.

With this diagram,
the crossing number of a Montesinos knot $K$ is given by the following.

%
% Theorem
%

\begin{theorem*}[Theorem 10 of \cite{LT}]
If a link $L$ admits an $n$-crossing, reduced Montesinos diagram,
then $L$ cannot be projected with fewer than $n$ crossings.
\end{theorem*}

In short, 
in order to calculate the crossing number of a given Montesinos knot $K$,
we have only to prepare a reduced Montesinos diagram
and count its crossings.

Thus, by putting them all together, we have the following.

%
% Proposition
%

\begin{proposition}
\label{Prop:RestrictedKnotsAndDiagrams}
For a Montesinos knot $K=K(R_1,R_2,\ldots,R_\NumTangles)$,
there is a tuple $(R^\prime_1,R^\prime_2,\ldots,R^\prime_\NumTangles)$
satisfying either $(A)$ or $(B)$,
such that the knot $K$ is isotopic to $K^\prime=K(R^\prime_1,R^\prime_2,\ldots,R^\prime_\NumTangles)$
and has a diagram $D=D(R^\prime_1,R^\prime_2,\ldots,R^\prime_\NumTangles)$.
Furthermore, the diagram $D$ is a reduced Montesinos diagram defined in \cite{LT} 
and attains the crossing number of $K$.
\end{proposition}

%%%%%%%%%%%%%%%%%%%%%%%%%%%%%%%%%%%%%%%%%%%
%
% Section : Proof of the main theorem 
%
%%%%%%%%%%%%%%%%%%%%%%%%%%%%%%%%%%%%%%%%%%%

\section{Proof of the main theorem}
\label{Sec:ProofOfTheMainTheorem}

The key of the proof is
the relation between 
the number of crossings of the alternating diagram of a rational tangle
and lengths of monotonic edgepaths for the rational tangle.
Roughly, summing up these relations gives the main theorem,
since the crossing number of a Montesinos knot
is calculated as the sum of the numbers of crossings of the alternating diagrams 
of rational tangles by \cite{LT}.

%
% Subsection
%

\subsection{Number of crossings and lengths of monotonic edgepaths}
\label{SubSec:CrossingsAndMonotonicEdgepaths}

We start with a single rational tangle.

Let $\Crossing(D)$ denote the number of crossings of a tangle diagram $D$.
The number $\Crossing(D(R))$ of crossings of the alternating tangle diagram $D(R)$ of a rational tangle $R$
is related to the lengths of monotonic edgepaths for $R$
by the following lemma.

Let $\BasicEdgepath_\inc$ and $\BasicEdgepath_\dec$ be
the monotonically increasing and decreasing basic edgepaths,
which are determined uniquely for given $R$.
By adding $\infty$-edges,
we obtain
the monotonically increasing and decreasing type III edgepaths
$\Edgepath_{\textrm{III},\inc}$ and $\Edgepath_{\textrm{III},\dec}$
from the basic edgepaths.
Note that,
since $\infty$-edges do not contribute to the twist,
we ignore $\infty$-edges when we judge a type III edgepath
is monotonic or not.
%
%%% Recall that
%%% we ignore $\infty$-edges when we judge 
%%% whether a type III edgepath system is monotonic or not.
%
According to the sign of $R$,
we can construct a monotonic type II edgepath 
with an ending vertex $\angleb{0}$,
which is
a monotonically increasing type II edgepath $\Edgepath_{\textrm{II},\inc}$ if $R<0$
and a monotonically decreasing type II edgepath $\Edgepath_{\textrm{II},\dec}$
if $R>0$.
%%% labeled $0$

Then, we obtain the following key lemma.

%
% Lemma
%

\begin{lemma}
\label{Lem:LengthsAndCrossing}
Assume that $R$ is a non-integral fraction.
Then,
%%%If $R>0$,
\begin{eqnarray*}
\Crossing(D(R))=|\Edgepath_{\textrm{III},\inc,\ge0}|+|\Edgepath_{\textrm{II},\dec}| \textrm{\ \ \ \ \ \ if $R>0$.}
%%%\label{Eq:IdentityForASingleTangleR>0}
\end{eqnarray*}
%%%If $R<0$,
\begin{eqnarray*}
\Crossing(D(R))=|\Edgepath_{\textrm{II},\inc}|+|\Edgepath_{\textrm{III},\dec,\ge0}| \textrm{\ \ \ \ \ \ if $R<0$.}
%%%\label{Eq:IdentityForASingleTangleR<0}
\end{eqnarray*}
Furthermore, 
\begin{eqnarray*}
\Crossing(D(R))=|\Edgepath_{\textrm{III},\inc,\ge0}|+|\Edgepath_{\textrm{III},\dec,\ge0}|
\textrm{\ \ \ \ \ \ if $|R|<1$.}
%%%\label{Eq:IdentityForASingleTangle|K|<1}
\end{eqnarray*}
$\Edgepath_{*,\ge 0}$ denotes the part of $\Edgepath_{*}$ lying in $u\ge 0$.
\end{lemma}

\begin{figure}[hbt]
 \begin{picture}(150,190)
%block-a-nor.eps PS 216x316 216x316+0+0 PseudoClass 256c 67kb 0.010u 0:01
 \put(0,0){
  \put(0,0){\scalebox{0.5}{\includegraphics{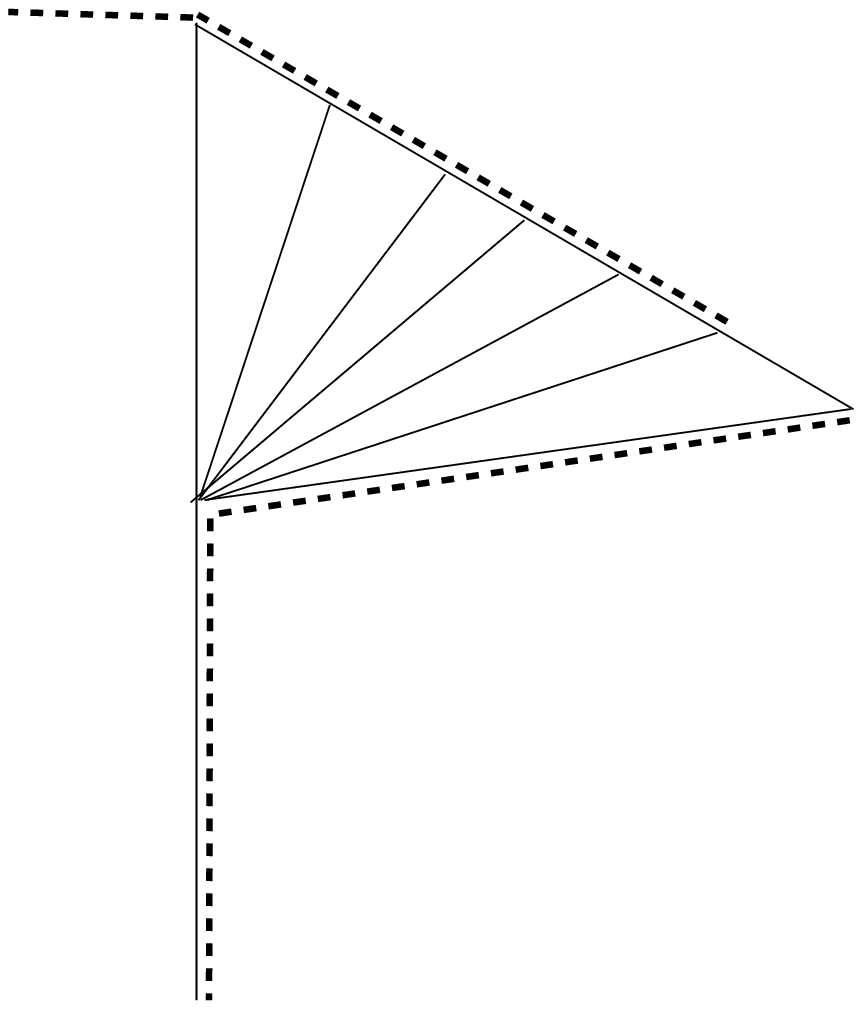}}}% 108x158
  %\put(72,149){\rotatebox{120}{\scalebox{0.5}{$\cdots$}}}

  \put(-13,177){$\angleb{\infty}$}
  \put(-3,173){\vector(1,-2){8}}
  \put(12,160){$\Edgepath_{\inc,1}$}

  \put(-12,12){$\angleb{(0)}$}
  \put(13,14){\vector(1,0){16}}

  \put(-17,85){$\angleb{(z_1)}$}
  \put(13,88){\vector(1,0){16}}

  \put(54,164){$\angleb{(z_1,1)}=\angleb{(z_1+1)}$}
  \put(52,165){\vector(-2,-1){15}}
  \put(76,135){$\Edgepath_{\inc,2}$}

  \put(71,151){$\angleb{(z_1,2)}$}
  \put(69,153){\vector(-2,-1){15}}

  \put(116,128){$\angleb{(z_1,z_2-1)}$}
  \put(113,129){\vector(-2,-1){15}}

  \put(128,116){$\angleb{(z_1,z_2)}$}
  \put(127,118){\vector(-2,-1){15}}

  \put(144,106){$\angleb{(z_1,z_2+1)}$}
  \put(143,106){\vector(-2,-1){15}}
  \put(78,83){$\Edgepath_{\dec,2}$}

  \put(37,50){$\Edgepath_{\dec,1}$}
 }

 \end{picture}
 \caption{$\Edgepath_{\inc,l}$ and $\Edgepath_{\dec,l}$ for $l=1$ and $l=2$}
 \label{Fig:EdgepathsInTriangleL=1or2}
\end{figure}

\begin{figure}[hbt]
 \begin{picture}(220,220)

%block-b-nor.eps PS 204x150 204x150+0+0 PseudoClass 256c 30kb 0.000u 0:01
 \put(0,120){
  \put(0,0){\scalebox{0.5}{\includegraphics{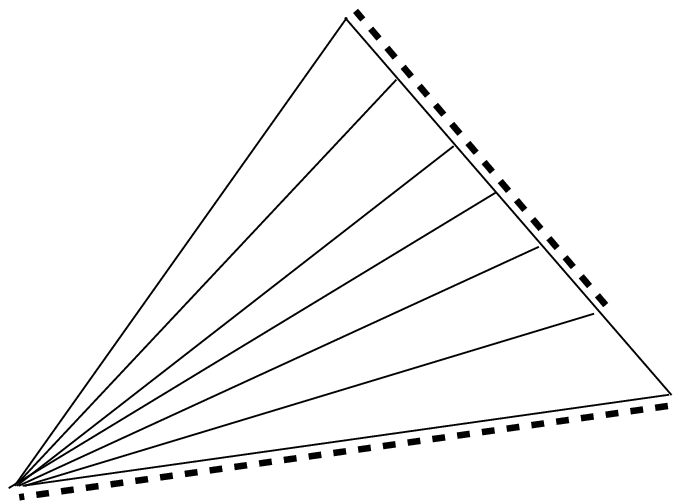}}}% 204x150
  %\put(72,9){\rotatebox{120}{\scalebox{0.5}{$\cdots$}}}

  \put(-78,24){$\angleb{(z_1,z_2,\ldots,z_{l-1})}$}
  \put(-5,21){\vector(1,-2){8}}

  \put(-35,84){$\angleb{(z_1,z_2,\ldots,z_{l-1},1)}=\angleb{(z_1,z_2,\ldots,z_{l-1}+1)}$}
  \put(71,79){\vector(-2,-1){15}}
  \put(73,53){$\Edgepath_{\inc,l}$}

  \put(79,69){$\angleb{(z_1,z_2,\ldots,z_{l-1},2)}$}
  \put(77,71){\vector(-2,-1){15}}

  \put(100,46){$\angleb{(z_1,z_2,\ldots,z_{l-1},z_l-1)}$}
  \put(98,45){\vector(-2,-1){15}}

  \put(109,34){$\angleb{(z_1,z_2,\ldots,z_{l-1},z_l)}$}
  \put(107,36){\vector(-2,-1){15}}

  \put(119,23){$\angleb{(z_1,z_2,\ldots,z_{l-1},z_l+1)}$}
  \put(116,24){\vector(-2,-1){15}}
  \put(50,0){$\Edgepath_{\dec,l}$}
 }

 \put(0,20){
  \put(0,0){\scalebox{0.5}{\includegraphics{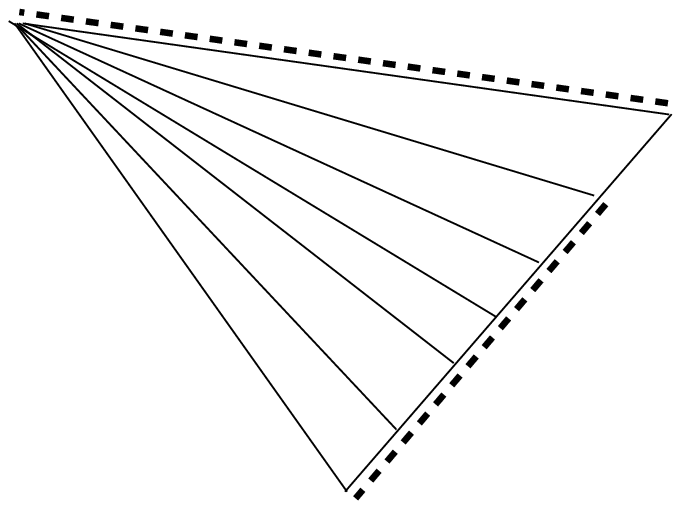}}}% 204x150
  %\put(72,71){\rotatebox{120}{\scalebox{0.5}{$\cdots$}}}

  \put(-78,46){$\angleb{(z_1,z_2,\ldots,z_{l-1})}$}
  \put(-5,54){\vector(1,2){8}}

  \put(-35,-13){$\angleb{(z_1,z_2,\ldots,z_{l-1},1)}=\angleb{(z_1,z_2,\ldots,z_{l-1}+1)}$}
  \put(71,-4){\vector(-2,1){15}}
  \put(75,21){$\Edgepath_{\dec,l}$}

  \put(79,1){$\angleb{(z_1,z_2,\ldots,z_{l-1},2)}$}
  \put(77,3){\vector(-2,1){15}}

  \put(100,27){$\angleb{(z_1,z_2,\ldots,z_{l-1},z_l-1)}$}
  \put(98,29){\vector(-2,1){15}}

  \put(108,39){$\angleb{(z_1,z_2,\ldots,z_{l-1},z_l)}$}
  \put(106,40){\vector(-2,1){15}}

  \put(118,51){$\angleb{(z_1,z_2,\ldots,z_{l-1},z_l+1)}$}
  \put(115,53){\vector(-2,1){15}}
  \put(50,73){$\Edgepath_{\inc,l}$}
 }
 \end{picture}
 \caption{$\Edgepath_{\inc,l}$ and $\Edgepath_{\dec,l}$ in a large triangle $T_l$
          for $3\le l \le k-1$}
 \label{Fig:EdgepathsInTriangleL=3tok-1}
\end{figure}

\begin{proof}
Assume that $R$ is positive, and 
$R$ is represented by a continued fraction $(z_1, z_2, \ldots, z_k)$
satisfying $k\ge 2$ and the condition (\ref{Eq:ConditionOfUniqueContinuedFraction}).
This continued fraction corresponds to the alternating diagram $D(R)$ of the rational tangle $R$.
We also prepare monotonic edgepaths $\Edgepath_{\textrm{III},\inc}$ and $\Edgepath_{\textrm{II},\dec}$.

\bigskip

%%%

We divide these monotonic edgepaths 
by use of triangles defined as follows.
For integers $2\le l\le k$ and $m\ge 1$,
there exists a large triangle with vertices
\[
\angleb{(z_1,z_2,\ldots,z_{l-1})}, \ 
\angleb{(z_1,z_2,\ldots,z_{l-1},1)}, \ 
\angleb{(z_1,z_2,\ldots,z_{l-1},m)}
.
\]
See Figure \ref{Fig:EdgepathsInTriangleL=3tok-1}.
Note here that an expression like $(z_1, z_2, \ldots, z_a)$ denotes a rational number $R$ corresponding to the continued fraction,
and hence,
$\angleb{(z_1, z_2, \ldots, z_a)}$ denotes a vertex $\angleb{R}$ of the diagram $\Diagram$.
Moreover, this triangle consists of some smaller triangles with vertices 
\[
\angleb{(z_1,z_2,\ldots,z_{l-1})}, \ %
\angleb{(z_1,z_2,\ldots,z_{l-1},j)},  \ %
\angleb{(z_1,z_2,\ldots,z_{l-1},j+1)}
\]
for $1\le j\le m-1$.
The existence of such a triangle is explained roughly as follows.

%%  First, we consider the set of vertices
%%  $\angleb{1}$, $\angleb{1/2}$, $\ldots$, $\angleb{1/t}$, $\ldots$, $\angleb{0}$.
%%  For a natural number $z$, it is easy to see that 
%%  there are the edges 
%%  $\angleb{0}$\,--\,$\angleb{1/z}$ and $\angleb{1/z}$\,--\,$\angleb{1/(z+1)}$.
%%  Vertices other than $\angleb{0}$ lie on the same line $u+v=1$.
%%  $1>1/2>\ldots>1/t>\ldots>0$ holds.
%%  Thus, a finite part of this gives a large triangle.
%%  %
%%  Next,
%%  for the set of vertices, or equivalently the set of fractions,
%%  we add the same integer to all fractions at the same time.
%%  With respect to natures mentioned as above,
%%  the operation preserves most of them.
%%  %
%%  Replacing each fraction by its inverse fraction at the same time
%%  reverses increasing-ness, decreasing-ness and magnitude relation.
%%  Though, most of other natures are preserved.
%%  %
%%  Repeating adding a fixed integer and replacing by inverse,
%%  we can confirm that such a large triangle above surely exists.
%
%
An important fact here is that
if $(z_1,z_2,\ldots,z_l)=p/q$ and $(z_1,z_2,\ldots,z_l,1)=r/s$ with $|ps-qr|=1$,
then $(z_1,z_2,\ldots,z_l,j)=\{(j-1)p+r\}/\{(j-1)q+s\}$ holds for any $j\ge 1$.
This is easily shown for $(0,j)$'s,
and then for $(z_{l},j)$'s, for $(0,z_{l},j)$'s, for $(z_{l-1},z_{l},j)$'s, and so on, inductively.
With the fact, the following features are easily confirmed.
There are edges $\angleb{(z_1,z_2,\ldots,z_l)}$\,--\,$\angleb{(z_1,z_2,\ldots,z_l,j)}$ and $\angleb{(z_1,z_2,\ldots,z_l,j)}$\,--\,$\angleb{(z_1,z_2,\ldots,z_l,j+1)}$.
For the sequence of vertices $\angleb{(z_1,z_2,\ldots,z_l,1)}$, $\angleb{(z_1$,$z_2$,$\ldots$,$z_l$,$2)}$, $\ldots$,
their $v$-coordinates are monotonic
and their $u$-coordinates are monotonically increasing.
For all edges of the form $\angleb{(z_1,z_2,\ldots,z_l,j)}$\,--\,$\angleb{(z_1,z_2,\ldots,z_l,j+1)}$, if we regard them as line segments,
then the slopes of the segments coincide with each other.
This is because
we have
\begin{eqnarray*}
 \angleb{(z_1, z_2, \ldots,z_l,j)}
 &=&\angleb{\frac{(j-1)p+r}{(j-1)q+s}},
 \\
 \angleb{(z_1, z_2, \ldots,z_l,j+1)}
 &=&\angleb{\frac{jp+r}{jq+s}},
\end{eqnarray*}
their $uv$-coordinates are
\[
(1-\frac{1}{(j-1)q+s},\frac{(j-1)p+r}{(j-1)q+s})
\text{~and~}
(1-\frac{1}{jq+s},\frac{jp+r}{jq+s}),
\]
and the slope of the segment between them is found to be
\[
 \frac{ps-qr}{q},
\]
which does not depend on $j$.
Hence, the vertices $\angleb{(z_1,z_2,\ldots,z_l,j)}$'s for $j\ge1$ lie on the same line.
Thus, we have confirmed the existence of the triangle mentioned as above.

%%%

\bigskip

Now, for $l=1$,
let $\Edgepath_{\inc,l}$$=$$\angleb{\infty}$\,--\,$\angleb{z_1+1}$,
$\Edgepath_{\dec,l}$$=$$\angleb{0}$\,--\,$\angleb{1}$\,--\,$\ldots$\,--\,$\angleb{z_1}$.
We have $|\Edgepath_{\dec,l}|=z_1$.

For $2\le l \le k-1$, 
we focus on the triangle $T_l$ having the vertices
\[
\angleb{(z_1,z_2,\ldots,z_{l-1})}, \ %
\angleb{(z_1,z_2,\ldots,z_{l-1},1)}, \ %
\angleb{(z_1,z_2,\ldots,z_{l-1},{z_l}+1)}
.
\]
See Figure \ref{Fig:EdgepathsInTriangleL=1or2} if $l=2$,
Figure \ref{Fig:EdgepathsInTriangleL=3tok-1} if $3\le l \le k-1$.
Which of the two large triangles
in Figure \ref{Fig:EdgepathsInTriangleL=3tok-1}
the triangle $T_l$ looks like
depends on the parity of $l$.
The large triangle $T_l$ includes two edgepaths
\begin{eqnarray*}
\Edgepath_{l,A}&=&
\angleb{(z_1,z_2,\ldots,z_{l-1},1)}\textrm{\,--\,}%
\angleb{(z_1,z_2,\ldots,z_{l-1},2)}\textrm{\,--\,}%
\\
&&\ldots\textrm{\,--\,}%
\angleb{(z_1,z_2,\ldots,z_{l-1},z_{l}-1)}\textrm{\,--\,}%
\angleb{(z_1,z_2,\ldots,z_{l-1},z_{l})},
\\
\Edgepath_{l,B1}&=&
\angleb{(z_1,z_2,\ldots,z_{l-1})}\textrm{\,--\,}%
\angleb{(z_1,z_2,\ldots,z_{l-1},z_{l}+1)}
%\\
%\Edgepath_{l,B2}&=&%
%\angleb{(z_1,z_2,\ldots,z_{l-1})}\textrm{\,--\,}%
%\angleb{(z_1,z_2,\ldots,z_{l-1},z_{l}+1)}
.
\end{eqnarray*}
One of the two edgepaths is monotonically increasing,
and the other is monotonically decreasing.
Note that $\Edgepath_{l,A}$ can be a point when $z_l =1$, 
and such an edgepath can be regarded as both increasing and decreasing.
Let $\Edgepath_{\inc,l}$ and $\Edgepath_{\dec,l}$ denote the increasing one and the decreasing one.
The important fact here is that
$|\Edgepath_{\inc,l}|+|\Edgepath_{\dec,l}|=z_l-1+1=z_l$.

\begin{figure}[hbt]
 \begin{picture}(220,220)
 \put(0,120){
  \put(0,0){\scalebox{0.5}{\includegraphics{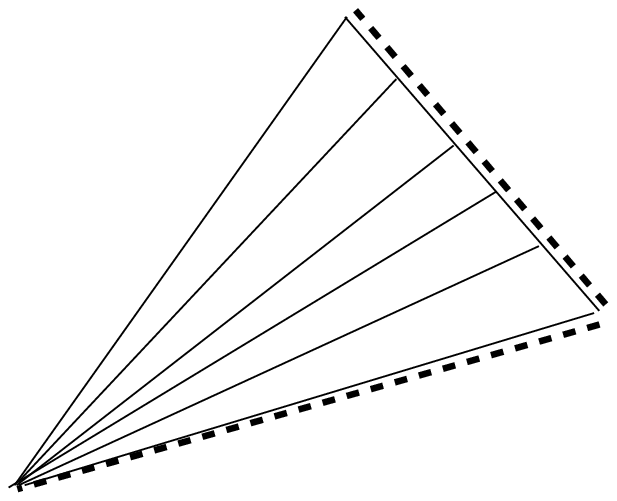}}}% 204x150
  %\put(72,9){\rotatebox{120}{\scalebox{0.5}{$\cdots$}}}

  \put(-78,24){$\angleb{(z_1,z_2,\ldots,z_{l-1})}$}
  \put(-5,21){\vector(1,-2){8}}

  \put(-35,84){$\angleb{(z_1,z_2,\ldots,z_{l-1},1)}=\angleb{(z_1,z_2,\ldots,z_{l-1}+1)}$}
  \put(71,79){\vector(-2,-1){15}}
  \put(73,53){$\Edgepath_{\inc,l}$}

  \put(79,69){$\angleb{(z_1,z_2,\ldots,z_{l-1},2)}$}
  \put(77,71){\vector(-2,-1){15}}

  \put(100,46){$\angleb{(z_1,z_2,\ldots,z_{l-1},z_l-1)}$}
  \put(98,45){\vector(-2,-1){15}}

  \put(109,34){$\angleb{(z_1,z_2,\ldots,z_{l-1},z_l)}$}
  \put(107,36){\vector(-2,-1){15}}
  \put(50,7){$\Edgepath_{\dec,l}$}

 }

 \put(0,20){
  \put(0,0){\scalebox{0.5}{\includegraphics{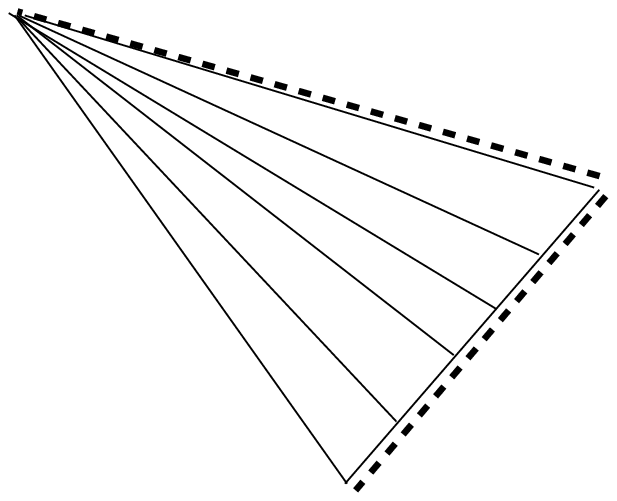}}}% 204x150
  %\put(72,141){\rotatebox{120}{\scalebox{0.5}{$\cdots$}}}

  \put(-78,46){$\angleb{(z_1,z_2,\ldots,z_{l-1})}$}
  \put(-5,54){\vector(1,2){8}}

  \put(-35,-13){$\angleb{(z_1,z_2,\ldots,z_{l-1},1)}=\angleb{(z_1,z_2,\ldots,z_{l-1}+1)}$}
  \put(71,-4){\vector(-2,1){15}}
  \put(75,21){$\Edgepath_{\dec,l}$}

  \put(79,1){$\angleb{(z_1,z_2,\ldots,z_{l-1},2)}$}
  \put(77,3){\vector(-2,1){15}}

  \put(100,27){$\angleb{(z_1,z_2,\ldots,z_{l-1},z_l-1)}$}
  \put(98,29){\vector(-2,1){15}}

  \put(108,39){$\angleb{(z_1,z_2,\ldots,z_{l-1},z_l)}$}
  \put(106,40){\vector(-2,1){15}}
  \put(50,66){$\Edgepath_{\inc,l}$}

 }

 \end{picture}
 \caption{$\Edgepath_{\inc,l}$ and $\Edgepath_{\dec,l}$ in a large triangle $T_l$
          for $l=k$}
 \label{Fig:EdgepathsInTriangleL=k}
\end{figure}

For $l=k$, 
we define 
$\angleb{(z_1,z_2,\ldots,z_{l-1})}$%
$\angleb{(z_1,z_2,\ldots,z_{l-1},1)}$%
$\angleb{(z_1,z_2,\ldots,z_{l-1},{z_l})}$
as the triangle $T_l$.
We use 
\begin{eqnarray*}
\Edgepath_{l,B2}&=&%
\angleb{(z_1,z_2,\ldots,z_{l-1})}\textrm{\,--\,}%
\angleb{(z_1,z_2,\ldots,z_{l-1},z_{l})}
\end{eqnarray*}
instead of $\Edgepath_{l,B1}$.
Though,
the situation is similar to the previous case,
and hence, the identity $|\Edgepath_{\inc,l}|+|\Edgepath_{\dec,l}|=z_l$ still holds.
See Figure \ref{Fig:EdgepathsInTriangleL=k}.

By paying attention to the fact that $(z_1,z_2,\ldots,z_l+1)=(z_1,z_2,\ldots,z_l,1)$,
we can find that $\Edgepath_{\textrm{III},\inc}$ is obtained by combining $\Edgepath_{\inc,l}$ $(1\le l \le k)$, and $\Edgepath_{\textrm{II},\dec}$ by combining $\Edgepath_{\dec,l}$ $(1 \le l \le k)$.
Furthermore, we have
\begin{eqnarray*}
|\Edgepath_{\textrm{III},\inc,\ge0}|+|\Edgepath_{\textrm{II},\dec}|
=\sum_{l=2}^{k} (|\Edgepath_{\inc,l}|+|\Edgepath_{\dec,l}|)+|\Edgepath_{\dec,1}|
=\sum_{l=1}^{k} z_l
=\Crossing(D(R))
.
\end{eqnarray*}
The argument is similar for negative $R$.

If $|R|<1$, then $z_1$ is $0$.
Hence, we have 
\begin{eqnarray*}
%% \label{Eq:IdentityForASingleTangle}
|\Edgepath_{\textrm{III},\inc,\ge 0}|+|\Edgepath_{\textrm{\textrm{III}},\dec,\ge 0}|=\Crossing(D(R)).
\end{eqnarray*}
\end{proof}

%
% Figure
%

%
% Figure
%

%
% Figure
%

\begin{figure}[hbt]
 \begin{picture}(130,250)
%47over36-inc-dec-2.eps PS 225x307 225x307+0+0 PseudoClass 256c 68kb 0.010u 0:01
 \put(0,0){
  \put(-21,0){\scalebox{0.8}{\includegraphics{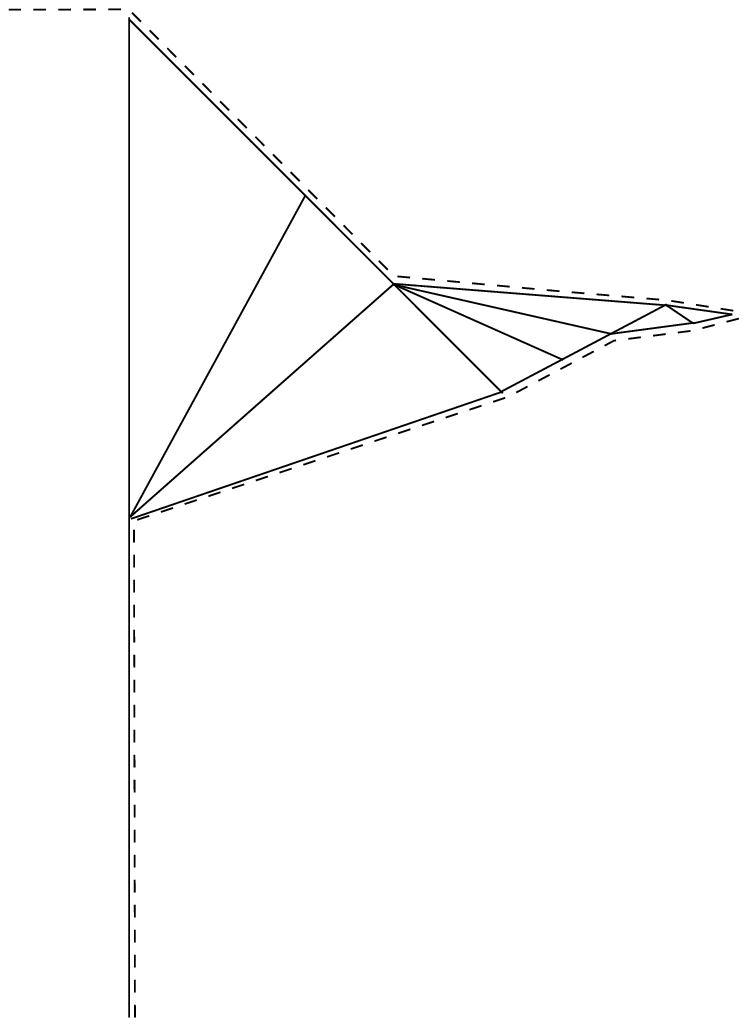}}} % 180x246

  \put(38,215){$\Edgepath_{\textrm{III},\inc}$}
  \put(38,123){$\Edgepath_{\textrm{II},\dec}$}
  \put(-10,0){$\angleb{0}$}
  \put(-10,120){$\angleb{1}$}
  \put(-30,247){$\angleb{\infty}$}
  \put(87,199){\vector(-1,-2){10}}
  \put(77,203){$\angleb{4/3}$}
  \put(86,126){\vector(1,2){10}}
  \put(72,115){$\angleb{5/4}$}

  \put(139,184){\vector(-1,-2){5}}
  \put(127,187){$\angleb{17/13}$}
  \put(123,150){\vector(0,1){10}}
  \put(106,139){$\angleb{13/10}$}

  \put(163,153){\vector(-2,1){20}}
  \put(148,142){$\angleb{30/23}$}

  \put(174,167){\vector(-1,0){20}}
  \put(176,165){$\angleb{47/36}$}
 }
 \end{picture}
 \caption{$\Edgepath_{\textrm{III},\inc}$ and $\Edgepath_{\textrm{II},\dec}$, for $R=47/36$.}
 \label{Fig:Edgepaths47over36}
\end{figure}

For instance,
in Figure \ref{Fig:Edgepaths47over36},
$\Edgepath_{\textrm{III},\inc}$ and $\Edgepath_{\textrm{II},\dec}$
are depicted for $47/36=(1,3,3,1,2)$.
In the figure,
for making it clear,
the horizontal coordinates of the vertices are somewhat altered from correct values,
though, it does not matter in the argument.
Note that $|\Edgepath_{\textrm{III},\inc,\ge0}|+|\Edgepath_{\textrm{II},\dec}|=4+6=10$ coincides with $\Crossing(D(47/36))=1+3+3+1+2=10$.

%
% Figure
%

By the argument in the above proof,
we obtained the following about the possible edgepaths.

%
% Proposition
%

\begin{proposition}
The twists of monotonic basic edgepath systems give the upper and lower bounds
for twists of basic edgepath systems, 
of type I edgepath systems and of type III edgepath systems.
\label{Prop:TwistOfMonotonicEdgepaths}
\end{proposition}

\begin{proof}
Let
$\BasicEdgepathSystem_\inc=(\BasicEdgepath_{\inc,1},\ldots,\BasicEdgepath_{\inc,\NumTangles})$ and
$\BasicEdgepathSystem_\dec=(\BasicEdgepath_{\dec,1},\ldots,\BasicEdgepath_{\dec,\NumTangles})$
denote
the unique monotonically increasing basic edgepath system
and the unique monotonically decreasing basic edgepath system.
$\EdgepathSystem=(\Edgepath_{1},\ldots,\Edgepath_{\NumTangles})$ denotes a basic edgepath system, a type I edgepath system or a type III edgepath system.
If
$\Twist(\BasicEdgepath_{\inc,i})\le 
 \Twist(\Edgepath_{i})\le
 \Twist(\BasicEdgepath_{\dec,i})$
holds for each $i$,
then
we immediately have
$\sum_{i=1}^\NumTangles \Twist(\BasicEdgepath_{\inc,i})\le
 \sum_{i=1}^\NumTangles \Twist(\Edgepath_{i})\le
 \sum_{i=1}^\NumTangles \Twist(\BasicEdgepath_{\dec,i})$,
or equivalently,
$\Twist(\BasicEdgepathSystem_\inc)\le 
 \Twist(\EdgepathSystem)\le
 \Twist(\BasicEdgepathSystem_\dec)$.
Hence,
it suffices to show that 
$\Twist(\BasicEdgepath_{\inc})\le 
 \Twist(\Edgepath)\le
 \Twist(\BasicEdgepath_{\dec})$
holds for each edgepath $\Edgepath=\Edgepath_i$.

For a fixed starting point $\angleb{R}=\angleb{p/q}$,
all non-constant edgepaths lie in the region $X$ 
bounded by
the monotonically increasing basic edgepath $\BasicEdgepath_\inc$,
the monotonically decreasing basic edgepath $\BasicEdgepath_\dec$,
and exactly one vertical edge
as illustrated in Figure \ref{Fig:Edgepaths47over36}.
In this figure,
$X$ is the union of $8$ triangles.
As seen from the figure,
roughly, the region $X$ consists of a sequence of triangles.
Hence, no vertices of $\Edgepath$ other than a rational endpoint
lie in the interior of $X$.
All increasing (complete) edges in the region end at a vertex on $\BasicEdgepath_\inc$. 
For each intermediate vertex on $\BasicEdgepath_\inc$,
there are one or more increasing edges ending at the vertex.
A single edgepath contains at most one of them, 
and just one of them is always included in $\BasicEdgepath_\inc$.
This implies that
the number of increasing edges in an edgepath 
is less than or equal to that of 
the monotonically increasing basic edgepath. 
At the same time,
the number of decreasing edges in an edgepath
is more than or equal to that of 
the monotonically increasing basic edgepath,
which is zero.
Recall that 
the twist can be calculated as
twice of 
the total length of decreasing non-$\infty$-edges
minus the total length of increasing non-$\infty$-edges.
Hence, we have $\Twist(\BasicEdgepath_\inc)\le\Twist(\Edgepath)$.
The argument is similar for the monotonically decreasing basic edgepath.
\end{proof}

% In short,
% the above lemma roughly tells that
% for a single rational tangle,
% the half of its contribution to the diameter of the set of boundary slopes
% coincides with the number of crossings of its alternating tangle diagram.

%
% Subsection
%

\subsection{Crossing number and diameter}

We sum up identities in Lemma \ref{Lem:LengthsAndCrossing}
to all rational tangles of a Montesinos knot,
and obtain the relation between the diameter and the crossing number.

We separate the argument according to which of the conditions (A) and (B) is satisfied.

%
% Subsubsection
%

\subsubsection{Alternating case}

If a tuple $(R_1,R_2,\ldots,R_\NumTangles)$ satisfies the condition (A),
then the Montesinos knot $K(R_1,R_2,\ldots,R_\NumTangles)$ is alternating.
Conversely, by the following well-known fact, 
together with the another well-known fact that any Montesinos knot without $1/0$-tangles is prime,
an alternating Montesinos knot is isotopic to a Montesinos knot 
with a tuple satisfying not (B) but (A).
%
%
% Proposition
%
\begin{proposition}[e.g., \cite{L}]
If a link $L$ has a connected, irreducible, alternating diagram with $n$ crossings,
then
the number of crossings of any diagram of the link $L$ is $n$ or greater.
Moreover, 
the number of crossings of any prime, non-alternating diagram is greater than $n$.
\end{proposition}

In this case, the situation is rather simple,
and we have the following lemma.

%
% Lemma
%

\begin{lemma}
\label{Lem:CrossingNumberAndDiameter:CommonSign}
For a tuple $(R_1,R_2,\ldots,R_\NumTangles)$,
assume that all $R_j$'s have the common sign.
Then, for the Montesinos knot $K=K(R_1,R_2,\ldots,R_\NumTangles)$, 
the following identity holds.
\begin{eqnarray*}
2\,\Crossing(K)&=&\Diam(K).
\label{Eq:Diam:Crossing:CommonSign}
\end{eqnarray*}
\end{lemma}

%
% Proof
%

\begin{proof}
%%% Note that such a Montesinos knot $K$ is alternating.
Assume here that all $R_j$'s are positive.
As mentioned above,
by Theorem 10 of \cite{LT},
$\Crossing(K)=\sum_{j=1}^{\NumTangles} \Crossing(D(R_j))$%
.
Let 
$\EdgepathSystem_{\textrm{III},\inc}$ 
and
$\EdgepathSystem_{\textrm{II},\dec}$
denote
the monotonically increasing type III edgepath system
and
the monotonically decreasing type II edgepath system
respectively
whose edgepaths are constructed in the same manner as in Subsection \ref{SubSec:CrossingsAndMonotonicEdgepaths}.
With respect to these,
by Lemma \ref{Lem:LengthsAndCrossing},
we have 
$\Crossing(D(R_j))=|\Edgepath_{\textrm{III},\inc,j,\ge 0}|+|\Edgepath_{\textrm{II},\dec,j}| \textrm{ for any $j$}$%
.
One complete edge contributes $\pm 2$ to the twist.
The twists of two edgepath systems $\EdgepathSystem_{\textrm{III},\inc}$ and $\EdgepathSystem_{\textrm{II},\dec}$ have opposite signs.
Combining these facts, we have:
\begin{eqnarray*}
&& 2\,\Crossing(K)=2\,\sum_{j=1}^{\NumTangles} \Crossing(D(R_j)) 
=2\,\sum_{j=1}^{\NumTangles} (|\Edgepath_{\textrm{III},\inc,j,\ge 0}|+|\Edgepath_{\textrm{II},\dec,j}|) \\
&&=\sum_{j=1}^{\NumTangles} 2\,|\Edgepath_{\textrm{III},\inc,j,\ge 0}|+\sum_{j=1}^{\NumTangles} 2\,|\Edgepath_{\textrm{II},\dec,j}| 
=|\Twist(\EdgepathSystem_{\textrm{III},\inc})|+|\Twist(\EdgepathSystem_{\textrm{II},\dec})| 
\\
&&=|\Twist(\EdgepathSystem_{\textrm{III},\inc})-\Twist(\EdgepathSystem_{\textrm{II},\dec})|.
\end{eqnarray*}

Now, both $\EdgepathSystem_{\textrm{III},\inc}$ and
$\EdgepathSystem_{\textrm{II},\dec}$ are minimal.
$\EdgepathSystem_{\textrm{III},\inc}$ gives an essential surface 
by Proposition 2.5 of \cite{HO}
since the sum of $v$-coordinates at $u=0$ of the edgepaths in the edgepath system $\EdgepathSystem_{\textrm{III},\inc}$ is at least $3$.
In addition,
$\EdgepathSystem_{\textrm{II},\dec}$ corresponds to an essential surface
by one of Proposition 2.4, 2.6, 2.7, 2.8 of \cite{HO}
since it is monotonic.

For a type I or type III edgepath system $\EdgepathSystem$,
we have 
$\Twist(\EdgepathSystem_{\textrm{III},\inc})
= \Twist(\BasicEdgepathSystem_\inc)
\le \Twist(\EdgepathSystem)
\le \Twist(\BasicEdgepathSystem_\dec)
\le \Twist(\EdgepathSystem_{\textrm{II},\dec})$
by Proposition \ref{Prop:TwistOfMonotonicEdgepaths}.
If $\EdgepathSystem=(\Edgepath_1,\Edgepath_2,\ldots,\Edgepath_\NumTangles)$
is type II,
there is another type II edgepath system 
$\EdgepathSystem^\prime=(\Edgepath^\prime_1,\Edgepath^\prime_2,\ldots,\Edgepath^\prime_\NumTangles)$
each of whose edgepath $\Edgepath^\prime_i$ shares a basic edgepath with $\Edgepath_i$
and ends at $\angleb{0}$.
Note that $\Twist(\EdgepathSystem)=\Twist(\EdgepathSystem^\prime)$ holds
since the edgepath systems share a basic edgepath system
and have the same total amount of signed length of vertical edges
by the gluing consistency.
By an argument similar to Proposition \ref{Prop:TwistOfMonotonicEdgepaths}
with considering effect of vertical edges on twist,
we have $\Twist(\EdgepathSystem_{\textrm{III},\inc})=\Twist(\BasicEdgepathSystem_\inc)\le\Twist(\EdgepathSystem^\prime)\le\Twist(\EdgepathSystem_{\textrm{II},\dec})$.
Hence,
$\Twist(\EdgepathSystem_{\textrm{III},\inc})\le\Twist(\EdgepathSystem)\le\Twist(\EdgepathSystem_{\textrm{II},\dec})$ holds for type II edgepath systems,
and hence for edgepath systems of any type.

In short,
$\EdgepathSystem_{\textrm{III},\inc}$
and
$\EdgepathSystem_{\textrm{II},\dec}$
are found to give the minimal and the maximal twists.
Finally,
\begin{eqnarray*}
2\,\Crossing(K)=|\Twist(\EdgepathSystem_{\textrm{III},\inc})-\Twist(\EdgepathSystem_{\textrm{II},\dec})|=\Diam(K).
\end{eqnarray*}
\end{proof}

%%%%%
%
% Subsubsection
%

\subsubsection{Non-alternating case}

Note that 
Lemma \ref{Lem:CrossingNumberAndDiameter:CommonSign} as it is cannot be extended to the general Montesinos knots. 
For instance, $K=K(1/2,1/3,-1/3)=8_{20}$ has non-integral diameter $38/3$.
If a tuple $(R_1,R_2,\ldots,R_\NumTangles)$ satisfies the condition (B),
then we have the following lemma.

%
% Lemma
%

\begin{lemma}
\label{Lem:CrossingNumberAndDiameter:DistinctSigns}
For a tuple $(R_1,R_2,\ldots,R_\NumTangles)$,
assume that 
$\{R_j\}$ includes both positive and negative fractions, and $|R_j|<1$ holds for any $j$.
Then, for the Montesinos knot $K=K(R_1,R_2,\ldots,R_\NumTangles)$,
the following inequality holds.

\begin{eqnarray*}
2\,\Crossing(K)&\ge&\Diam(K)
.
\label{Eq:Diam:Crossing:DistinctSigns}
\end{eqnarray*}
\end{lemma}

\begin{proof}
First, we prepare monotonic type III edgepath systems $\EdgepathSystem_{\textrm{III},\inc}$ and $\EdgepathSystem_{\textrm{III},\dec}$.
Similarly to the argument in the proof of the previous lemma,
we have
\begin{eqnarray*}
2\,\Crossing(K)&=&|\Twist(\EdgepathSystem_{\textrm{III},\inc})-\Twist(\EdgepathSystem_{\textrm{III},\dec})|.
\end{eqnarray*}
Note that, the condition $|R_j|<1$ is required to apply Lemma \ref{Lem:LengthsAndCrossing}.

We think about the range of twists.
The twists of the monotonic type III edgepath systems are
$\Twist(\EdgepathSystem_{\textrm{III},\inc})
=\Twist(\BasicEdgepathSystem_{\inc})$
and
$\Twist(\EdgepathSystem_{\textrm{III},\dec})
=\Twist(\BasicEdgepathSystem_{\dec})$.
Though, the surfaces corresponding to these edgepath systems may be inessential.
We check that
$\Twist(\EdgepathSystem_{\textrm{III},\inc})$
and $\Twist(\EdgepathSystem_{\textrm{III},\dec})$
are the lower and upper bounds of the twist $\Twist$.
The claim is true for type I and type III edgepath systems by Proposition \ref{Prop:TwistOfMonotonicEdgepaths}.
For a type II edgepath system $\EdgepathSystem$,
by preparing $\EdgepathSystem^\prime$ as in the proof of the previous lemma,
we easily have 
$\Twist(\EdgepathSystem)=\Twist(\EdgepathSystem^\prime)$
and
$\Twist(\EdgepathSystem_{\textrm{III},\inc})
\le\Twist(\EdgepathSystem^\prime)
\le\Twist(\EdgepathSystem_{\textrm{III},\dec})$
by an argument similar to Proposition \ref{Prop:TwistOfMonotonicEdgepaths}.
Hence, the twist is bounded by
$\Twist(\EdgepathSystem_{\textrm{III},\inc})$ and $\Twist(\EdgepathSystem_{\textrm{III},\dec})$
also for any type II edgepath system.
Thus,
\begin{eqnarray*}
|\Twist(\EdgepathSystem_{\textrm{III},\inc})-\Twist(\EdgepathSystem_{\textrm{III},\dec})|
&\ge& \Diam(K)
.
\end{eqnarray*}

Then, we have
\begin{eqnarray*}
2\,\Crossing(K)&\ge &\Diam(K)
.
\end{eqnarray*}
\end{proof}

%
% Theorem
%

Finally, combining
two-bridge cases in \cite{MMR},
Proposition \ref{Prop:RestrictedKnotsAndDiagrams},
Lemma \ref{Lem:CrossingNumberAndDiameter:CommonSign}
and Lemma \ref{Lem:CrossingNumberAndDiameter:DistinctSigns},
we have Theorem \ref{Thm:CrossingNumberAndDiameter:Main}.

%%%%%%%%%%%%%%%%%%%%%%%%%%%%%%%%%%%%%%%%%%%
%
% References
%
%%%%%%%%%%%%%%%%%%%%%%%%%%%%%%%%%%%%%%%%%%%

%%%%%%%%%%%%%%%%%%%%%%%%%%%%%%%%%%%%%%%%%%%
%
% Appendix , Trash
%
%%%%%%%%%%%%%%%%%%%%%%%%%%%%%%%%%%%%%%%%%%%

%%% \appendix

\end{document}